\documentclass[11pt,a4paper]{article}
\usepackage[utf8]{inputenc}
\usepackage{epic,eepic,epsfig,amsmath,amssymb,amsthm}
\usepackage{t1enc}

\usepackage[francais,english]{babel}

\usepackage{caption}

\usepackage{color}

\usepackage{bbm}

\def\virgp{\raise 2pt\hbox{,}}

\renewcommand{\geq}{\geqslant}
\renewcommand{\leq}{\leqslant}
\def\N{{\mathbb N}}

\def\R{{\mathbb R}}

\def\virgp{\raise 2pt\hbox{,}}
\def\cdotpv{\raise 2pt\hbox{;}}

\def\1{\mathbbm{1}}

\newtheorem{theorem}{Theorem}[section]
\newtheorem{corollary}[theorem]{Corollary}
\newtheorem{proposition}[theorem]{Proposition}
\newtheorem{lemma}[theorem]{Lemma}

\newtheorem{pte}[theorem]{Property}

\theoremstyle{remark}
\newtheorem{remark}{Remark}[section]

\theoremstyle{definition}
\newtheorem{definition}{Definition}[section]

\newtheorem*{notation}{Notation}

\theoremstyle{definition}

\theoremstyle{definition}

\begin{document}
	
	\title{The Finite Volume Method on Sierpi\'{n}ski Simplices}
	
	\author{Nizar Riane$^\dag$, Claire David$^\ddag$}

	\maketitle
	\centerline{$^\dag$ Universit\'e Mohammed V de Rabat, Maroc\footnote{nizar.riane@gmail.com} }
	\vskip 0.5cm
	
	\centerline{$^\ddag$ Sorbonne Universit\'e}
	
	\centerline{CNRS, UMR 7598, Laboratoire Jacques-Louis Lions, 4, place Jussieu 75005, Paris, France\footnote{Claire.David@Sorbonne-Universite.fr}}

\begin{abstract}
	In the sequel, we extend the Strichartz average approach for the Laplacian, which uses average values of a function over basic sets, rather than pointwise ones as classically done in the literature, and goes back to the seminal work of~S.~Kusuoka and~X.~Y.~Zhou, to Sierpi\'{n}ski Simplices. To our knowledge, the finite volume method had not been applied before to fractals themselves, with their intrinsic Laplacian. We build the finite volume scheme for the heat equation with Dirichlet boundary conditions, give its consistency error, a stability condition of Courant-Friedrichs-Lewy type, obtained by means of an inclusion between the spectrum of the finite volume matrix and the one of the graph Laplacian, and prove the convergence of the scheme. Numerical results are given in the case of the Sierpi\'{n}ski Gasket and Tetrahedron, and compared with the finite difference method.
\end{abstract}

	\maketitle
	\vskip 1cm
	

	\vskip 1cm
	
	\vskip 1cm
	
	\vskip 1cm
	
\section{Introduction}

\hskip 0.5cm Why should one focus on fractals as, for instance, Sierpi\'{n}ski simplices? Let us go back to the study of Rammal Rammal, the condensed matter physicist~\cite{RammalToulouse1983}, who laid the emphasis upon the fact that fractals could ``fill the gap between well-ordered crystalline structures and those that are disordered''. Then, the work of~S.~Havlin and D.~Ben~Avraham~\cite{HavlinAvraham1987} showed that diffusion in those disordered systems did not follow classical laws, and thus required specifically designed tools.\\

In the years that followed, people began to take an interest in the applications of fractals, for numerous reasons. To begin with, one may note that fractal structures, because of the very large size of the length or surface they cover, show insulating or repellent properties. One may refer to silica aerogel, with a branched microstructure of the fractal type, which is a remarkable insulator, in so far as it almost completely stops the propagation of heat. In chemistry, to go on, scientists started to design dendritic polymers with an arborescent structure (as it is the case of arborol), where patterns such as the Sierpi\'{n}ski Gasket, the Koch curve, or the Menger sponge, are respectively used as models for linear polymer chains, or porous media~\cite{Liu1986},~\cite{HavlinAvraham1987}. In recent years, particularly, one has seen the design of polymers of fractal type (see, for instance,~\cite{FractalPolymersSG}, or~\cite{ChimieSierpinski}, where the authors describe a terpyridine-based architecture which simulates
a first-generation Sierpi\'nski triangle, as it can be observed on Figures~1 and~2).\\

Dendritic polymers are at stake in diffusion phenomena; for instance, they are now more and more used in pharmacology, since their very specific architecture enables one to control (prevent or facilitate) drug delivery.

\begin{figure}[h!]
	\center{\psfig{height=6cm,width=14cm,angle=0,file=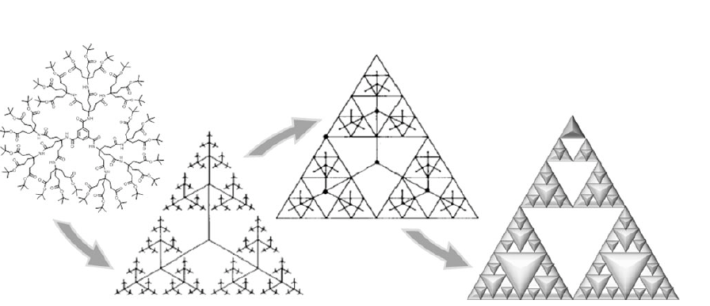}}\\
	\caption{Conceptual progression of a $1 \to 3$  dendritic branching pattern
		and its geometric relationship to the classical Sierpi\'nski triangle~\cite{ChimieSierpinski}. Figure reprinted with permission from the authors and publisher.}
\end{figure}

\begin{figure}[h!]
	\center{\psfig{height=8cm,width=10cm,angle=0,file=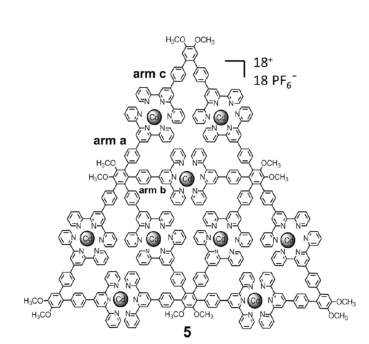}}\\
	\caption{Terpyridine-based Sierpi\'nski triangle~\cite{ChimieSierpinski}. Figure reprinted with permission from the authors and publisher.}
\end{figure}

One thus requires specifically fitted numerical tools. In the case of the Sierpi\'{n}ski gasket~$\mathfrak{SG}$, K.~Dalrymple, R.~S.~Strichartz, and J.~Vinson~\cite{DalrympleFDM} studied the heat and wave equations by means of the spectral decomposition of the solution (heat kernel), which involves eigenvalues and eigenvectors, and also introduced the corresponding finite difference approximation on the graphs approximating~$\mathfrak{SG}$. It was followed by the work of M.~Gibbons, A.~Raj and R.~S.~Strichartz~\cite{GibbonsFEM}, who gave a description of how one may build approximate solutions, by means of piecewise harmonic, or biharmonic, splines, i.e. the finite element method, again in the case of~$\mathfrak{SG}$. In~our work~\cite{RianeDavidFDMS}, following the seminal work of~J.~Kigami and~R.~S.~Strichartz for Laplacians on Sierpi\'nski simplices, we built and analyzed the related finite difference scheme on Sierpi\'nski simplices of any dimension.\\

As for us, our contribution laid in defining the discretization of the considered PDE's (heat and wave equation), by taking into account the recursive construction of the matrix related to a sequence of graph Laplacians. We thus did not call for approximations of the eigenvalues. This enabled us not only to compute the consistency error, but, also, to set stability conditions of Courant-Friedrichs-Lewy type, and, then, to prove the convergence of the scheme. \\

The extension to the finite volume method (FVM) of those numerical methods specifically devoted to fractal differential operators was a natural step, especially in the light of the paper of~R.~S.~Strichartz~\cite{Strichartz2001}, where he shows how the symmetric Laplacian on~$\mathfrak{SG}$ can be defined entirely in terms of average values of a function over basic sets (see also~\cite{RavierStrichartz2015}), rather than using pointwise ones as classically done in the literature. The use of average values goes back to the seminal work of~S.~Kusuoka and~X.~Y.~Zhou~\cite{KusuokaZhou1992} for the Sierpi\'nski Carpet. This approach is interesting, in so far as the Sierpi\'nski Carpet, contrary to Sierpi\'nski simplices, is not what may be called a post-critically finite fractal, i.e. a set which can be disconnected by removing a finite number of points, and therefore requires a non-pointwise valued approach, as in the case of the FVM method. The advantages of dealing with the FVM method consist, first, of its physical approach, which adapts to any geometry, and, in our case, to a fractal one, while being, at the same time, a conservative approach, and the basis of several general numerical codes. To our knowledge, the finite volume method had not been applied before to fractals themselves, with their intrinsic Laplacian, and, in particular, to Sierpi\'{n}ski simplices.\\

We hereafter present our results, in the case of Sierpi\'{n}ski simplices (Gasket and Tetrahedron), for the heat equation with Dirichlet boundary conditions: we give, first, the scheme and its numerical implementation, then, estimates of the consistency error, a stability condition of Courant-Friedrichs-Lewy type, obtained by means of an inclusion between the spectrum of the finite volume matrix and the one of the graph Laplacian, the convergence of the scheme, and a comparison with the finite-difference method.

	\section{Framework of the study}

	\begin{notation} 
		
		\noindent We will denote by~$\N$ the set of natural integers,
		\noindent and set:
		
		$$\N^\star= \N\setminus \left \lbrace 0 \right \rbrace \, \cdot$$

	\end{notation}
	
	\vskip 1cm
	\hskip 0.5cm	In the following, given an integer~\mbox{$d\geq 2$}, we place ourselves in the Euclidean space of dimension~$d-1$, referred to a direct orthonormal frame. The usual Cartesian coordinates will be denoted by~$(x_1,x_2,\hdots ,x_{d-1})$.\\

	\vskip 1cm
	
	\subsection{Graph approximation}
	
	\begin{notation}  
		In the following,~\mbox{$\left \lbrace f_1,\hdots, f_d \right \rbrace $} is a set of contractive maps, of ratio~$\displaystyle \frac{1}{2}$, where, for any integer~$i$ of~\mbox{$\left \lbrace 1, \hdots, d\right \rbrace $},~\mbox{$P_i\, \in\,\R^{d-1} $} is the fixed point of~$f_i$, the points~$P_1,\hdots,P_d$ being the vertices of a regular simplex. For any~$X$ of~$\R^{d-1}$:

		$$f_i(X)=\frac{1}{2}(X+P_{i})$$
	\end{notation}

	\vskip 1cm
	
	\begin{theorem}\textbf{(Hutchinson~\cite{Hutchinson1981}, Barnsley and Demko~\cite{Barnsley1985})}\\

		\noindent Given a complete metric space~$(E,\delta)$, and a set~\mbox{$\left \lbrace f_i\right \rbrace_{1 \leq i \leq d }$}  of contractions on~$E$ with respect to the metric~$\delta$, there exists a unique non-empty compact subset~$K\subset E$ such that:
		
	\begin{equation}
	 	K =  \underset{  i=1}{\overset{d }{\bigcup}}\,f_{i}\left ( {K} \right ) \, \cdot
	 \end{equation}

		
		\noindent The set~$K $ is said~\textbf{self-similar} with respect to the family~\mbox{$\left \lbrace f_1, \hdots, f_{d } \right \rbrace$}, and called attractor of the iterated function system (IFS)~\mbox{$\left \lbrace f_1 , \hdots, f_{d }  \right \rbrace$}.

	\end{theorem} 
	
	\vskip 1cm

	\begin{corollary}
		
		\noindent There exists a unique non-empty compact subset~$\mathfrak{SS} \subset \R^{d-1}$ such that:
			\begin{equation}\mathfrak{SS} = \underset{  i=1}{\overset{d}{\bigcup}}\, f_i(\mathfrak{SS})
				\end{equation}
		\noindent which will be called the Sierpi\'{n}ski Simplex.\\

	\end{corollary}

	\vskip 1cm

	\begin{definition}\textbf{Boundary (or initial) graph}\\
		
		\noindent We will denote by~$V_0$ the ordered set of the (boundary) points $\left \lbrace P_{1},\hdots ,P_{d  }\right \rbrace$, where $P_i$ is the fixed point of the contraction~$f_i$. The set~$V_0$, where, for any~$i$ of~\mbox{$\left \lbrace  1,...,d  \right \rbrace$}, the point~$P_i$ is adjacent to the other points~$P_{j}$,~$j\neq i$, constitutes a {complete graph}, that we will denote by~$  {{{\mathfrak {SS}}}_0}$.\\
		
		\noindent $V_0$ is called the {set of vertices} of the graph~$ {{{\mathfrak {SS}}}_0}$.
		
	\end{definition}

	\vskip 1cm

	\begin{definition}\textbf{$m^{th}$ order graph,~$m\,\in\, \N^\star$}\\

		\noindent For any strictly positive integer~$m$, we set:
			\begin{equation}V_m =\underset{  i=1}{\overset{d }{\bigcup}}\, f_i \left (V_{m-1}\right )\, \cdot  
			\end{equation}

		\noindent The set of points~$V_m$, where the points of the $m^{th}$-order cells~$f_{\cal W}(V_0)$ (see definition~\ref{Word} below) are linked in the same way as in~$\mathfrak{SS}_0$, is a (non oriented) graph, which we will denote by~$  {\mathfrak  {SS}}_m$. $V_m$ is called the {set of vertices} of the graph~$  {\mathfrak  {SS}}_m$. 
		
		\noindent By extension, we will write:

			\begin{equation}{\mathfrak  {SS}}_m=\underset{  i=1}{\overset{d}{\bigcup}}\, f_i \left ({\mathfrak  {SS}}_{m-1}\right )\, \cdot  
		\end{equation}
		
	\end{definition}

	\vskip 1cm
	
	\begin{pte}
		
		\vskip 0.5cm
		\noindent For any natural integer~$m$:
		\begin{equation}V_m \subset V_{m+1} \, \cdot  
		\end{equation}
		
	\end{pte}

	\vskip 1cm

	\begin{definition}\textbf{Word}\\
		\label{Word}
		
		\noindent Given a strictly positive integer~$m  $, we will call \textbf{number-letter} any integer~${\cal W}_i$ of~\mbox{$\left \lbrace 1, \hdots, d   \right \rbrace $}, and \textbf{word of length~$|{\cal W}|=m$}, on the graph~$ {\mathfrak {SS}}_m $, any sequence of number-letters of the form:

			\begin{equation}{\cal W}=\left ( {\cal W}_1, \hdots, {\cal W}_m\right) \, \cdot  
		\end{equation}
		
		\noindent We set:
		
		\begin{equation}f_{\cal W}= f_{{\cal W}_1} \circ \hdots \circ  f_{{\cal W}_m} \, \cdot  
	\end{equation}
		
		\noindent The sets~$f_{\cal W}(\mathfrak{SS})$, for~$|{\cal W}|=m$, are called the~$m^{th}$-order cells of~$\mathfrak{SS}$.
		
	\end{definition}
	
	\vskip 1cm
	
	\begin{definition}\textbf{Vertex }\\
		
		\noindent A point~$X$ of~${\mathfrak {SS}}$ will be called \textbf{\emph{vertex}} of the graph~$ {\mathfrak {SS}}$ if there exists a natural integer~$m$ such that:
		
		\begin{equation}  X \,\in\, V_m \, \cdot  
		\end{equation}
		
	\end{definition}
	
	\vskip 1cm
	
	\ 
		
		\begin{definition}\textbf{Consecutive (or neighbor) vertices}\\
			
			\noindent Two points~$X$ and~$Y$ of~${\mathfrak {SS}}$ will be called \textbf{\emph{consecutive vertices}}, or  \textbf{\emph{neighbor vertices}}, if there exist a natural integer~$m$, a word~${\cal W}$ of length~$m$, and two distinct integers~$i$ and~$j$ of~\mbox{$\left \lbrace  1,\hdots,d \right \rbrace$}, such that:
			
				\begin{equation}  X = f_{\cal W}(P_i) \quad \text{and} \quad Y = f_{\cal W}(P_j) \, \cdot  
			\end{equation}
			
			\noindent In other words,~$X$ and~$Y$ are two distinct vertices of the same~$m^{th}$-order cell.
			
		\end{definition}

	\vskip 1cm
	
	\begin{definition}\textbf{Edge relation}\\
		
		\label{Edge}
		
		\noindent For any~\mbox{$m\,\,\in\,\N$}, two points~$X$ and~$Y$ of~$ {\mathfrak {SS}}_m$ will be called \emph{\textbf{adjacent}} if and only if there exist a word~{${\cal W}=\left ( {\cal W}_1, \hdots, {\cal W}_m\right)$} of length~$ m$, and two distinct integers~$i$ and~$j$ of~\mbox{$\left \lbrace  1,\hdots,d \right \rbrace$}, such that~$X=f_{\cal W}(P_i)$ and~$Y=f_{\cal W}(P_j)$, i.e. if and only if~$X$ and~$Y$ both belong to the iterate:

	\begin{equation}f_{\cal W} \,V_0=\left (f_{{\cal W}_1} \circ \hdots \circ  f_{{\cal W}_m} \right) \,V_0 \, \cdot  
		\end{equation}
		
		\noindent We set:
		
		\begin{equation} X \underset{m }{\sim}  Y \, \cdot  
		\end{equation}

		\noindent Given two points~$X$ and~$Y$ of~$ {\mathfrak {SS}} $, we will say that~$X$ and~$Y$ are \textbf{\emph{adjacent}} if and only if there exists a natural integer~$m$ such that:
		
		\begin{equation}X  \underset{m }{\sim}  Y \, \cdot  
	\end{equation}
	\end{definition}
	
	\vskip 1cm

	\begin{definition}\textbf{Addresses}\\
		
		\noindent For any~\mbox{$m\,\,\in\,\N$}, and any vertex~$X$ of~${\mathfrak {SS}}_m$, we will call {address of the vertex~$X$} an expression of the form

		\begin{equation}X= f_{ \cal W} \left (P_i\right ) 
		\end{equation}

		\noindent	where~$\cal W$ is a word of length~$m$, and~$i$ a natural integer in~\mbox{$\left \lbrace 1,\hdots, d \right \rbrace $}.
	\end{definition}
	
	\begin{remark}
		The address of a vertex is not unique: any point~$X$ of~$V_m\setminus V_0$ is the junction point of exactly two~$m^{th}$-order cells, and thus admits two addresses~$X=f_{\cal W}(P_i)=f_{{\cal W}'}(P_j)$ of length~$m$. As a consequence, any point~$X$ of~$V_m\setminus V_0$ has exactly~$2\,(d-1)$ neighbours in~$\mathfrak{SS}_m$.
	\end{remark}

	\vskip 1cm

	\begin{proposition}
		\noindent For any~\mbox{$m\,\,\in\,\N$}, we will denote by~$\mathcal{N}_m$ the number of vertices of the graph ${\mathfrak {SS}}_m$. We have:

	\begin{equation} \mathcal{N}_0  =d\, \quad , \quad \forall\, m\,\,\in\,\N^\star\, : \quad 
		\mathcal{N}_m   =d\,\mathcal{N}_{m-1}- \frac{d\,(d-1)}{2}=\displaystyle\frac{d^{m+1}+d}{2} \, \cdot  
		\end{equation}
	\end{proposition}

	\vskip 1cm
	
	\begin{proof}
		The graph $\mathfrak{SS}_m$ is the union of~$d$ copies of the graph $\mathfrak{SS}_{m-1}$. Each copy shares exactly one vertex with each one of the other copies. So, one may consider the copies as the vertices of a complete graph~$K_d$, the number of edges of which is equal to~$\displaystyle\frac{d\,(d-1)}{2}$, which leads to~$\displaystyle\frac{d\,(d-1)}{2}$ vertices counted twice. The closed form follows by induction.
	\end{proof}
	
	\vskip 1cm

	\begin{definition}

		\noindent For any~\mbox{$m\,\,\in\,\N$}, we consider the graph~$\mathcal{SS}_m$, built from $\mathfrak{SS}_m$ in the following way:
		\begin{enumerate}
			\item[\emph{i}.] each~$m^{th}$-order cell~$f_{\cal W}(\mathfrak{SS})$,~$|{\cal W}|=m$, becomes a vertex of~$\mathcal{SS}_m$ ;
			
			\item[\emph{ii}.] two vertices of~$\mathcal{SS}_m$ are linked if the corresponding cells have a vertex in common.
		\end{enumerate}
		
	\end{definition}
	
	\begin{pte}
		The graph~$\mathcal{SS}_m$ has~$d^m$ vertices. Each of its vertices has at most~$d$ neighbours: a cell~$f_{\cal W}(\mathfrak{SS})$ has~$d$ vertices, each one of them being shared with exactly one other cell, except the points of~$V_0$, which belong to a single cell. The~$d$ cells which contain a point of~$V_0$ thus have exactly~$d-1$ neighbours.
	\end{pte}
	
	\vskip 1cm

	\begin{remark}{\textbf{For the 2-Simplex (Triangle):\\}}
		\begin{figure}[!htb] 
			
			\minipage{0.4\textwidth}
			\includegraphics[width=\linewidth]{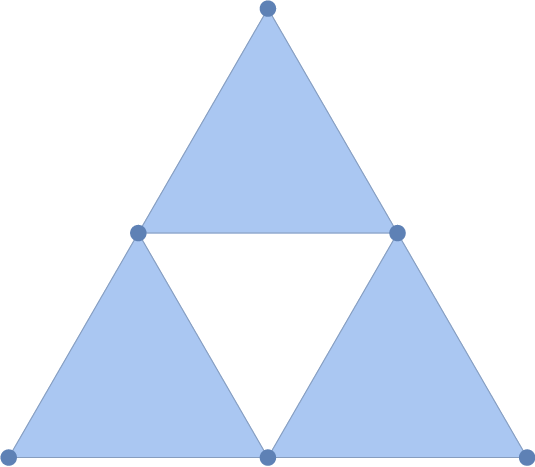}
			\caption{$\mathfrak{SS}_1$}
			\endminipage\hfill
			\minipage{0.4\textwidth}
			\includegraphics[width=\linewidth]{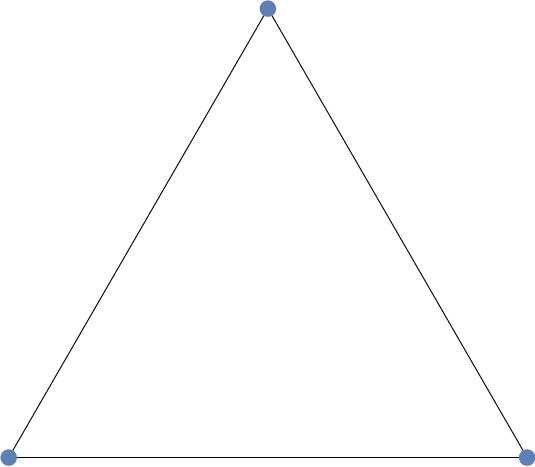}
			\caption{$\mathcal{SS}_1$}
			\endminipage\hfill
		\end{figure}
		
		\begin{figure}[!htb]
			
			\minipage{0.4\textwidth}
			\includegraphics[width=\linewidth]{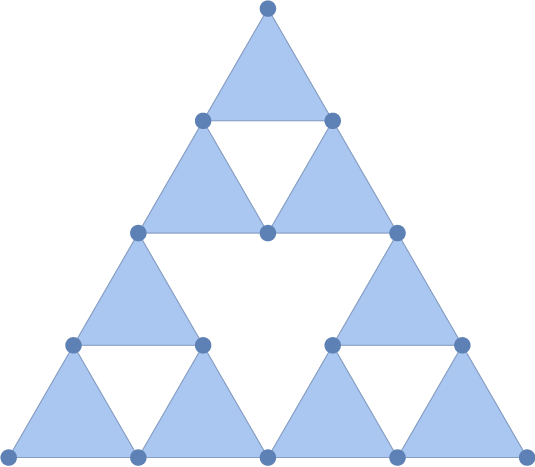}
			\caption{$\mathfrak{SS}_2$}
			\endminipage\hfill
			\minipage{0.4\textwidth}
			\includegraphics[width=\linewidth]{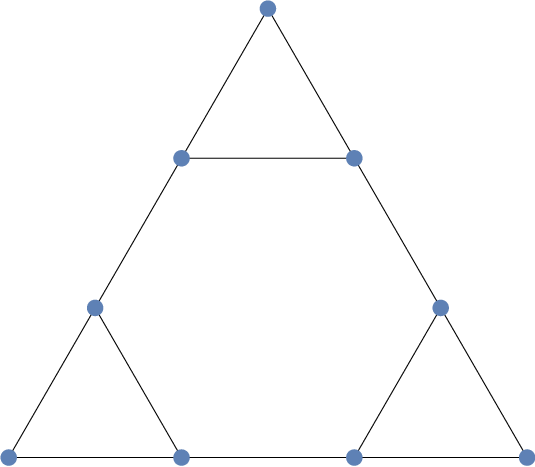}
			\caption{$\mathcal{SS}_2$}
			\endminipage\hfill
		\end{figure}
		
	\end{remark}
	
	\vskip 1cm
	
	\subsection{Self-similar measures on Sierpi\'nski Simplices}
	
	\begin{definition}\textbf{Self-similar measure, on the Sierpi\'nski simplex~$\mathfrak {SS}$}\\
		
		\noindent A Borel probability measure~$\mu$ on~$\R^{d-1}$, supported by~$\mathfrak {SS}$, will be said to be \textbf{self-similar} on the Sierpi\'nski simplex~$\mathfrak {SS}$, if there exists a family of strictly positive weights~\mbox{$\left (\mu_i\right)_{1 \leq i \leq d}$} such that:
		
		\begin{equation} \mu= \displaystyle \sum_{i=1}^{d} \mu_i\,\mu\circ f_i^{-1} \quad, \quad \displaystyle \sum_{i=1}^{d} \mu_i =1 \, \cdot  
		\end{equation}
		
		\noindent Given such a family of weights, there exists a unique self-similar probability measure satisfying this relation. For further precisions on self-similar measures, we refer to the works of~J.~E.~Hutchinson~(see \cite{Hutchinson1981}).
		
	\end{definition}
	
	\vskip 1cm
	
	\begin{pte}\textbf{Building of a self-similar measure, for the Sierpi\'nski simplex~$\mathfrak {SS}$}\\
		
		\noindent The Dirichlet forms introduced in the sequel require a positive Radon measure with full support.
		
		\noindent Let us denote by~$R_i=\displaystyle\frac{1}{2}$,~$1\leq i\leq d$, the contraction ratios of the similarities~$f_i$. We set, for any integer~$i$ belonging to~$\left \lbrace 1, \hdots, d \right \rbrace$:
		
		\begin{equation}\mu_i=R_i^{D_{H}\left({ \mathfrak {SS} }\right)}  
		\end{equation}
		
		\noindent where 
		
		$$D_{H}\left( { \mathfrak {SS} }\right)=\frac{\ln d}{\ln 2}$$
		
		\noindent is the Hausdorff dimension of the Sierpi\'nski simplex~$\mathfrak {SS}$. Since the open set condition is satisfied, the Hausdorff dimension is the solution of the Moran equation:
		
		\begin{equation} \displaystyle \sum_{i=1}^d R_i^{D_{H}\left(\mathfrak {SS}\right)}=1 \, \cdot  
	\end{equation}
		
		\noindent One may thus define a self-similar measure~$\mu$ on~$\mathfrak {SS}$ through:
		\begin{equation}\mu =\displaystyle \,\sum_{i=1}^d \mu_i  \mu \circ f_i^{-1}  
		\end{equation}
		
		\noindent Since~$\mu_i=\left(\frac{1}{2}\right)^{\frac{\ln d}{\ln 2}}=\displaystyle\frac{1}{d}$, this simply yields the standard measure:
		
			\begin{equation} \mu= \displaystyle \frac{1}{d}\displaystyle \sum_{i=1}^{d} \mu\circ f_i^{-1} \, \cdot  
			\end{equation}
		
		\noindent In particular, any~$m^{th}$-order cell has measure:
		
		$$\forall\,{\cal W}\,\in\,\{1,\hdots,d\}^m\,:\quad \mu\left(f_{\cal W}(\mathfrak{SS})\right)=d^{-m}$$
	\end{pte}

	\vskip 1cm

\subsection{Laplacians, on Sierpi\'nski Simplices  (we refer to~\cite{StrichartzLivre2006})}

	\begin{definition}{\textbf{Energy, on the graph~${{\mathfrak {SS}} }_m$,~$m \,\in\,\N$, of a pair of functions}}\\

		\noindent For any~\mbox{$m\,\,\in\,\N$}, and two real valued functions~$u$ and~$v$, defined on the set~$V_m$ of the~vertices of~${{\mathfrak {SS}} }_m$, we introduce \textbf{the energy, on the graph~${{\mathfrak {SS}} }_m$, of the pair of functions~$(u,v)$}, as:



	\begin{equation} {\cal{E}}_{{{\mathfrak {SS}} }_m}(u,v)= \displaystyle \sum_{X  \underset{m }{\sim}  Y}\left (u(X)-u(Y)\right )\,\left(v(X)-v(Y)\right)
		\end{equation}

		\noindent where the adjacency relation~$\underset{m}{\sim}$ has been defined in definition~\ref{Edge}, and where the sum is taken over the edges of~${{\mathfrak {SS}} }_m$, each edge being counted once.

	\end{definition}

	\vskip 1cm














	\begin{definition}{Dirichlet form, on a finite set (see~\cite{Kigami2003})}\\

		\noindent Let~$V$ denote a finite set. A \emph{\textbf{Dirichlet form}} on~$V$ is a symmetric bilinear form~${\cal E}$, defined on the space of real-valued functions on~$V$, such that:\\

		\begin{enumerate}

			\item For any real valued function~$u$ defined on~$V$:  ${\cal E}(u,u) \geq 0$.

			\item   $  {\cal {E}}(u,u)= 0$ if and only if~$u$ is constant on~$V$.

			\item For any real-valued function~$u$ defined on~$V$, if:~$ u_\star = \min\, (\max(u, 0) , 1)  $,
			i.e.:

			\begin{equation}\forall \,p \,\in\,V \, : \quad u_\star(p)= \left \lbrace \begin{array}{ccc} 1 & \text{if}& u(p) \geq 1 \\u(p) & \text{if}& 0 <u(p) < 1 \\0  & \text{if}& u(p) \leq 0 \end{array} \right. \,   \end{equation}

			\noindent then: ${ \cal{E}}(u_\star,u_\star)\leq { \cal{E}}(u,u)$ (Markov property).

		\end{enumerate}

	\end{definition}

	\vskip 1cm

	\begin{pte}

		\vskip 0.5cm

		 \noindent Given a natural integer~$m$, the symmetric bilinear form which, to any pair~$(u,v)$ of real-valued functions defined on~$V_m$, associates:

			\begin{equation} {\cal{E}}_{ {{\mathfrak {SS}}}_m}(u,v)= \displaystyle \sum_{X  \underset{m }{\sim}  Y}\left (u(X)-u(Y)\right )\, \left (v(X)-v(Y)\right )    \end{equation}

		\noindent is a Dirichlet form on~$ {{\mathfrak {SS}} }_m $. In particular, since the graph~${{\mathfrak {SS}} }_m$ is connected:

		\begin{equation}{\cal{E}}_{{{\mathfrak {SS}} }_m }(u,u)=0 \Leftrightarrow u  \text{ is constant} \, \cdot   \end{equation}

	\end{pte}

	\vskip 1cm
	\begin{proposition}

		\noindent For any strictly positive integer~$m$, if~$u$ is a real-valued function defined on~$V_{m-1}$, its \textbf{harmonic extension}, denoted by~$ \tilde{u}$, is obtained as the extension of~$u$ to~$V_m$ which minimizes the energy:

		\begin{equation}  {\cal{E}}_{{{\mathfrak {SS}}}_m}(\tilde{u},\tilde{u})=\displaystyle \sum_{X \underset{m }{\sim} Y} \displaystyle  (\tilde{u}(X)-\tilde{u}(Y))^2 \, \cdot   \end{equation}

	\end{proposition}

	\vskip 1cm

	\begin{remark} \textbf{More explicitly:}\\

		\noindent The link between~$   {\cal{E}}_{{{\mathfrak {SS}} }_m}$ and~$  {\cal{E}}_{ {{\mathfrak {SS}}}_{m-1}}$ is obtained through the introduction of two strictly positive constants~$r_m$ and~$r_{m-1}$ such that:

	\begin{equation}  r_{m }\, \displaystyle \sum_{X \underset{m  }{\sim} Y} \displaystyle  (\tilde{u}(X)-\tilde{u}(Y))^2  =   r_{m-1} \,\displaystyle  \sum_{X \underset{m-1 }{\sim} Y} (u(X)-u(Y))^2  \, \cdot   \end{equation}

		\noindent In particular:

		\begin{equation}  r_{1 }\, \displaystyle \sum_{X \underset{1  }{\sim} Y} (\tilde{u}(X)-\tilde{u}(Y))^2 =  r_{0}\,\displaystyle  \sum_{X \underset{0 }{\sim} Y} (u(X)-u(Y))^2 \, \cdot   \end{equation}

		\noindent We set:~$r_0=1$. Thus:

	\begin{equation}  {\cal{E}}_{ {{\mathfrak {SS}}}_1 }(\tilde{u},\tilde{u})= \displaystyle \frac{1}{ r_{1 } }\,  {\cal{E}}_{{{\mathfrak {SS}}}_0}(u,u) \, \cdot   \end{equation}

		\noindent Let us introduce:

	\begin{equation} r = \displaystyle \frac{1}{r_{1 }}   \end{equation}

		\noindent In the case of Sierpi\'nski simplices, a direct computation of the harmonic extension from~$V_0$ to~$V_1$ yields (see~\cite{StrichartzLivre2006}):

		\begin{equation} r=\displaystyle\frac{d}{d+2}\,\in\,\left]0,1\right[ \, \cdot \end{equation}

		\noindent and:

	\begin{equation}   {\cal{E}}_{m}(u)=r_m\, \sum_{X \underset{m }{\sim} Y}  (\tilde{u}(X)-\tilde{u}(Y))^2 \, \cdot   \end{equation}

		\noindent Since the determination of the harmonic extension of a function appears to be a local problem, on the graph~${{\mathfrak {SS}} }_{m-1}$, which is linked to the graph~${{\mathfrak {SS}} }_m$ by a similar process as the one that links~${{\mathfrak {SS}} }_1$ to~${{\mathfrak {SS}} }_0$, one deduces, for any strictly positive integer~$m$:

	\begin{equation} {\cal{E}}_{ {{\mathfrak {SS}}}_m }(\tilde{u},\tilde{u})= \displaystyle \frac{1}{ r_{1 }}\,  {\cal{E}}_{ {{\mathfrak {SS}}}_{m-1}}(u,u)\, \cdot   \end{equation}

		\noindent By induction, one gets:

	\begin{equation}  r_m=r_1^m =\displaystyle \frac{1}{r^m}\, \cdot   \end{equation}

		\noindent If~$v$ is a real-valued function, defined on~$V_{m-1}$, of harmonic extension~$ \tilde{v}$, we set:

	\begin{equation} {\cal{E}}_{m}(u,v)=\displaystyle \frac{1}{r^m} \, \sum_{X \underset{m }{\sim} Y}
		(\tilde{u}(X)-\tilde{u}(Y)) \, (\tilde{v}(X)-\tilde{v}(Y)) \, \cdot   \end{equation}

		\noindent For further details on the construction and existence of harmonic extensions, we refer to~\emph{\cite{Sabot1987}}.

	\end{remark}

	\vskip 1cm

	\begin{definition}\textbf{Renormalized energy, for a function~$u$, defined on~$V_m$,~$m\,\in\,\N$}\\

		\noindent Given a natural integer~$m$, one defines the \textbf{renormalized energy}, for a real-valued function~$u$, defined on~$V_m$, by:

		\begin{equation}
		{\cal{E}}_m (u)=    \displaystyle \sum_{X  \underset{m }{\sim}  Y} \displaystyle \frac{1}{r^m} \,
		\left (u (X)-u(Y)\right )^2 \, \cdot   \end{equation}

	\end{definition}

	\vskip 1cm

	\begin{definition}\textbf{Normalized energy, for a function~$u$, defined on~$V_\star$} \\

		Given a function~$u$ defined on~\mbox{$\displaystyle {V_\star =\underset{{i\in \N}}\bigcup \,V_i}$}, the sequence~$\left({\cal{E}}_m (u_{\mid V_m})\right)_{m\in\N}$ is non-decreasing. One defines the \textbf{normalized energy}:

	\begin{equation}
		{\cal{E}} (u)=
		\displaystyle \lim_{m \to + \infty} {\cal{E}}_m (u_{\mid V_m})\,\in\,[0,+\infty] \, \cdot
   \end{equation}

	\end{definition}
	\vskip 1cm

	\begin{definition}\textbf{Dirichlet form, for a pair of continuous functions defined on~$ {{\mathfrak {SS}} }  $} \\

		\noindent We define the Dirichlet form~$\cal{E}$ which, to any pair of real-valued, continuous functions~$(u,v)$ defined on~$  {{\mathfrak {SS}} }   $, of finite energy, associates:

	\begin{equation}
		{\cal{E}} (u,v)=
		\displaystyle \lim_{m \to + \infty}\displaystyle \sum_{X  \underset{m }{\sim}  Y} \displaystyle \frac{1}{r^m} \,
		\left (u_{\mid V_m}(X)-u_{\mid V_m}(Y)\right )\,\left(v_{\mid V_m}(X)-v_{\mid V_m}(Y)\right) \, \cdot   \end{equation}

		\noindent This limit exists, by polarization.

	\end{definition}

	\vskip 1cm

	\begin{notation}

		\noindent  We will denote by:

		\begin{enumerate}
			\item[\emph{i}.]  ~$\text{dom}\,{\cal E}$ the subspace of continuous functions defined on~$ {{\mathfrak {SS}}}$, such that:~\mbox{$\mathcal{E}(u)< + \infty$}.

			\item[\emph{ii}.] ~$\text{dom}_0\,{ \cal E}$ the  subspace of continuous functions defined on~$ {{\mathfrak {SS}}}$, which take the value zero on~$V_0$, and such that:~\mbox{$\mathcal{E}(u)< + \infty$}.

		\end{enumerate}

	\end{notation}
	\vskip 1cm

	\begin{lemma}

		The map:

	\begin{equation}\begin{array}{ccc}
		{\text{dom} \,{ \cal E}}  \big / {\text{Constants}}  \times {\text{dom} \,{ \cal E}}  \big / {\text{Constants}} &  \to& \R  \\
		(u,v) & \mapsto &
		{\cal{E}} (u,v) \end{array}  \end{equation}

		\normalsize
		\noindent defines an inner product on~\mbox{$	{\text{dom} \,{ \cal E}}  \big / {\text{Constants}} $}.

	\end{lemma}

	\vskip 1cm

	\begin{theorem}{\cite{StrichartzLivre2006}\\}

		$\left ( {\text{dom} \,{ \cal E}}  \big / {\text{Constants}} , { \cal E} ( \cdot, \cdot ) \right) $ is a Hilbert space.

	\end{theorem}

	\vskip 1cm

	\begin{definition}\textbf{Graph Laplacian of order~$m\,\in\,\N^\star$}\\

		\noindent For any strictly positive integer~$m$, and any real-valued function~$u$, defined on the set~$V_m$ of the vertices of the graph~${\mathfrak {SS}}_m$, we introduce the {graph Laplacian of order~$m$},~$\Delta_m(u)$, by:

	\begin{equation}\forall\, X\,\in\, V_m\setminus V_0 \, \quad \Delta_m u(X) =   \displaystyle\sum_{Y \in V_m,\,Y\underset{m}{\sim} X}
		\left (u(Y)-u(X)\right) \, \cdot \end{equation}

	\end{definition}

	\vskip 1cm

	\begin{definition}\textbf{Harmonic functions}\\

		\noindent A real-valued function~$u$, defined on~$\displaystyle {V_\star =\underset{{i\in \N}}\bigcup \,V_i}$, will be said to be \textbf{harmonic} if, for any natural integer~$m$, its restriction~$u_{\vert V_m}$ is harmonic:

	 	\begin{equation}\forall\, m\,\in\,\N^\star,\,\forall\, X\,\in\, V_m\setminus V_0 \, :\quad \Delta_m u_{\vert V_m}(X) =0 \, \cdot \end{equation}

		\noindent Equivalently, for any~$m\,\in\,\N^\star$,~$u_{\vert V_m}$ is the harmonic extension of~$u_{\vert V_{m-1}}$. A harmonic function is thus entirely determined by its values on~$V_0$.

	\end{definition}

	\vskip 1cm

	\begin{definition}

		\label{Lapl}
		\noindent We will denote by~$\text{dom}\, \Delta$ the existence domain of the Laplacian, on~$ {{\mathfrak {SS}}}$, as the set of functions~$u$ of~$\text{dom}\, \mathcal{E}$ such that there exists a continuous function on~$ {{\mathfrak {SS}}}$, denoted by~$\Delta \,u$, that we will call \textbf{Laplacian of~$u$}, such that,~\mbox{$\text{for any } v \,\in \,\text{dom}_0\, \mathcal{E}$}:

		\begin{equation}\begin{array}{ccc}\mathcal{E}(u,v)& =&
		\displaystyle \lim_{m \to + \infty}\displaystyle \sum_{X  \underset{m }{\sim}  Y} r^{-m}\,
		\left (u_{\mid V_m}(X)-u_{\mid V_m}(Y)\right )\,\left(v_{\mid V_m}(X)-v_{\mid V_m}(Y)\right)\\
		&=&-\displaystyle \int_{{\mathfrak {SS}}}  v\, \Delta u   \,d\mu \end{array}  \, \cdot \end{equation}

	\end{definition}

	\vskip 1cm

	\begin{theorem}{\cite{StrichartzLivre2006}\\}

		\begin{equation}\left (u \, \in\,  \text{dom}\, \Delta \quad \text{and} \quad \Delta \, u=0\right)
		\quad \text{\underline{if and only if}~$u$ is harmonic} \, \cdot \end{equation}

	\end{theorem}

	\vskip 1cm

	\begin{notation}
		\vskip 0.5cm

		\emph{i}.~${\cal H}_0\subset \text{dom}\, \Delta$ denotes the space of harmonic functions, i.e. the space of functions~$u \,\in\,\ \text{dom}\, \Delta$ such that:~\mbox{$\Delta\,u=0$}.\\

		\emph{ii}. Given a natural integer~$m$,~${\cal S} \left ({\cal H}_0,V_m \right)$ denotes the space of spline functions "of level~$m$",~$u$, defined on~$ {{\mathfrak {SS}}}$, continuous, such that, for any word~$\cal W$ of length~$m$,~\mbox{$u \circ f_{\cal W}$} is harmonic (see Definition~1.4.3 of~\cite{StrichartzLivre2006}), i.e.:

			\begin{equation} u \circ f_{\cal W}\,\in\,{\cal H}_0 \, \cdot \end{equation}

		\noindent A spline function of level~$m$ is thus entirely determined by its values on~$V_m$: it is obtained by prescribing these values, and by extending them harmonically to~$V_{m'}$, for any~$m'>m$.

	\end{notation}

	\vskip 1cm

	\begin{pte}

		For any natural integer~$m$:~${\cal S} \left ({\cal H}_0,V_m \right )\subset  \text{dom }{ \cal E}$.

	\end{pte}

	\vskip 1cm

	\begin{theorem}\textbf{Pointwise formula (Theorem~2.2.1 of~\cite{StrichartzLivre2006})}\\

		Let~$X \,\in \,V_\star \setminus V_0$, and, for any integer~$m$ such that~$X\,\in\,V_m$, let~\mbox{$\psi_X^{m}\,\in\,{\cal S} \left ({\cal H}_0,V_m \right)$} be the spline function such that:

			\begin{equation}\forall\, Y\,\in \,V_m\,:\quad \psi_X^{m}(Y)=\delta_{XY}  \quad,  \quad \text{where} \quad    \delta_{XY} =\left \lbrace \begin{array}{ccc}1& \text{if} & X=Y\\ 0& \text{else} &  \end{array} \right. \, \cdot \end{equation}

		\noindent The function~$\psi_X^{m}$ vanishes on any~$m^{th}$-order cell which does not contain~$X$, since it is there the harmonic extension of zero values. On each of the two~$m^{th}$-order cells containing~$X$, it is harmonic, with values in~$]0,1[$ outside the vertices of the cell. Its support is the union of the two~$m^{th}$-order cells containing~$X$, and one has:

		$$\int_{{\mathfrak {SS}}} \psi_X^{m}\,d \mu=\frac{2}{d^{m+1}}\quad , \quad r^{-m}\, \left \lbrace \int_{{\mathfrak {SS}}} \psi_X^{m}\,d \mu\right \rbrace^{-1}=\frac{d}{2}\,(d+2)^m \, \cdot$$

		\begin{enumerate}
			\item[ {i}.] For any function~$u$ of~$\text{dom}\, \Delta $, the sequence

			$$\left ( r^{-m}\, \left \lbrace \int_{{\mathfrak {SS}}} \psi_X^{m}\,d \mu\right \rbrace^{-1} \,\Delta_m u(X) \right)_{m\in\N}$$

			converges towards~$\Delta  u(X)$, uniformly with respect to~$X\,\in\,V_\star\setminus V_0$.

			\item[ {ii}.]  Conversely, given a continuous function~$u$ on~${{\mathfrak {SS}}}$ such that the sequence

			$$\left ( r^{-m}\, \left \lbrace \int_{{\mathfrak {SS}}} \psi_X^{m}\,d \mu\right \rbrace^{-1} \,\Delta_m u(X) \right)_{m\in\N}$$

			converges uniformly on~$V_\star \setminus V_0$ towards a function which extends continuously to~${{\mathfrak {SS}}}$, one has:

				\begin{equation}  u\, \in\, \text{dom}\, \Delta \quad \text{and} \quad
			\Delta \, u(X)=\displaystyle \lim_{m \to + \infty}  r^{-m}\, \left \lbrace \int_{{\mathfrak {SS}}} \psi_X^{m}\,d \mu\right \rbrace^{-1} \,\Delta_m u(X)   \, \cdot \end{equation}

		\end{enumerate}

	\end{theorem}

	\vskip 1cm

	\begin{definition}\textbf{Normal derivative (see Definition~2.3.1 of~\cite{StrichartzLivre2006})}\\

		\noindent Given a continuous function~$u$ on~${\mathfrak {SS}}$, and a boundary point~$P_i\,\in\,V_0$,~\mbox{$1\leq i \leq d$}, we will say that the normal derivative~$\partial_n u(P_i)$ exists if the limit
		$$ \partial_n u(P_i)=\displaystyle \lim_{m\rightarrow +\infty} \displaystyle \frac{1}{r^m}\, \sum_{Y \underset{m}{\sim} P_i} \left(u(P_i)-u(Y)\right)$$
		\noindent exists.\\

		\noindent More generally, given a boundary point~$X=f_{\cal W} (P_i)$ of a cell~$f_{\cal W}\left ( { \mathfrak {SS}} \right )$,~\mbox{${\cal W}\, \in\, \{1,\hdots ,d\}^\ell$}, the \textbf{local normal derivative} of~$u$ at~$X$, relative to the cell~$f_{\cal W}\left ( { \mathfrak {SS}} \right )$, is defined, when it exists, by:

		$$ \partial_n u(X)=\displaystyle \lim_{m\rightarrow +\infty} \displaystyle \frac{1}{r^m}\, \sum_{\substack{Y \underset{m}{\sim} X \\ Y\in f_{\cal W} \left ( { \mathfrak {SS}} \right )}} \left(u(X)-u(Y)\right)$$

		\noindent and satisfies:

		\begin{equation} \partial_n u(X) = r^{- \ell } \, \partial_n (u\circ f_{\cal W} ) (P_i)  \, \cdot \end{equation}

	\end{definition}

	\vskip 1cm

	\begin{theorem}\textbf{Gauss-Green formula (Theorem~2.3.2 of~\cite{StrichartzLivre2006})}\\

		\noindent Given~\mbox{$u\, \in \, \text{dom}\,\Delta$},~$\partial_n u(X)$ exists for all~\mbox{$X\, \in \,V_0$}, and:

		\begin{equation}  \mathcal{E}(u,v)=-\displaystyle \int_{\mathfrak {SS}} v\,\Delta u \,d\mu+ \displaystyle \sum_{X\,\in\,V_0} v(X)\,\partial_n u(X)    \end{equation}

		\noindent holds for all~\mbox{$v\, \in\,\text{dom} \, \mathcal{E}$}.
	\end{theorem}

	\vskip 1cm

	\begin{corollary}\textbf{(Theorem~2.4.1 of~\cite{StrichartzLivre2006})}\\

		\noindent Given~$u$ and~$v$ in~\mbox{$\text{dom}\,\Delta$}:

		\begin{equation}\displaystyle  \int_{\mathfrak {SS}} v\,\Delta u \, d\mu -\displaystyle \int_{\mathfrak {SS}} u\,\Delta v\, d\mu=\displaystyle  \sum_{X\,\in\,V_0} \left( v(X)\,\partial_n u(X) -u(X) \,\partial_n v(X) \right )  \end{equation}

	\end{corollary}

	\vskip 1cm

	\begin{theorem}\textbf{Matching condition (Theorem~2.3.3 of~\cite{StrichartzLivre2006})}\\

		\noindent Given~\mbox{$u\, \in \, \text{dom}\,\Delta$}, at each junction point

		$$X=f_{\cal W}(P_i)=f_{{\cal W} '}(P_j)\quad , \quad
		(i,j)\,\in\,\{1,\hdots,d\}^2 \quad , \quad
		\left ( {\cal W} , {\cal W}'\right) \, \in\, \{1,\hdots,d\}^m \times \{1,\hdots,d\}^m \quad , \quad {\cal W}\neq{\cal W}'$$
		\noindent the local normal derivatives relative to the two cells~$f_{\cal W}(\mathfrak{SS})$ and~$f_{{\cal W}'}(\mathfrak{SS})$ exist, and:
		\begin{equation} \partial_n u(f_{{\cal W}}(P_i))+\partial_n u(f_{{\cal W}'}(P_j))=0  \, \cdot \end{equation}

		\end{theorem}

	\vskip 1cm
	
	\begin{theorem}\textbf{Gluing (Theorem~2.5.1 of~\cite{StrichartzLivre2006})}\\
		
		\noindent Let~$m\,\in\,\N^\star$. Let us consider, for any word~${\cal W}$ of length~$m$, continuous functions~$u_{\cal W}$ and~$g_{\cal W}$, defined on the cell~$f_{\cal W}(\mathfrak{SS})$, such that~$\Delta u_{\cal W}=g_{\cal W}$ on~$f_{\cal W}(\mathfrak{SS})$, and let~$u$ and~$g$ be the functions obtained by gluing the pieces~$\{u_{\cal W}\}$ and~$\{g_{\cal W}\}$. Then:
		
		$$u\,\in\,\text{dom}\,\Delta \quad \text{and} \quad \Delta u=g$$
		
		\noindent if and only if the pieces~$\{u_{\cal W}\}$ and~$\{g_{\cal W}\}$ take the same values at each junction point (i.e.~$u$ and~$g$ are continuous), and the matching conditions on the normal derivatives hold at each junction point:
		
		$$\partial_n u_{\cal W}(f_{{\cal W}}(P_i))+\partial_n u_{{\cal W}'}(f_{{\cal W}'}(P_j))=0 \quad \text{for}\quad f_{{\cal W}}(P_i)=f_{{\cal W}'}(P_j)\,,\quad {\cal W}\neq{\cal W}' \, \cdot$$
		
	\end{theorem}

	\vskip 1cm

	\newpage
	
\section{The Finite Volume Method}

	\begin{notation} In the following,~$T$ is a strictly positive real number, while~$N$ is a non-zero natural integer. Let us introduce:

		\begin{equation}\displaystyle{h =\frac{T}{N}}\quad , \quad t_n=n\times h\quad, \quad  n=0,1,\hdots ,N
		\, \cdot  \end{equation}

	\end{notation}

	\subsection{The heat equation}

	\subsubsection{Formulation of the problem}

	\hskip 0.5cm We hereafter consider a solution~$u$ of the problem:

	\begin{equation}
	\label{eq:HeatFVM}
	\left \lbrace  \begin{array}{ccccc}
	\displaystyle \frac{\partial u}{\partial t}(t,X)-\Delta u(t,X)&=&0 & \forall (t,X)\,\in \,\left]0,T\right[  \times {\mathfrak{SS}}\\
	u(t,X)&=&0 & \forall \, ( t,X) \,\in \,\left[0,T\right[ \times \partial {\mathfrak{SS} } \\
	u(0,X)&=&g(X)  & \forall \,X\,\in \,{\mathfrak{SS}}
	\end{array}\right. 	\, \cdot  \end{equation}

	\noindent where~$\partial {\mathfrak{SS} }=V_0$.


	\vskip 1cm

	\begin{definition}
		\noindent Let~$m\, \in \,\N^\star$. We introduce the~$m^{th}$-\textbf{control volumes} as the~$d^m$ $m^{th}$-order cells
			\begin{equation} C^j_m=f_{{\cal W}^j}(\mathfrak{SS})\quad ,  \quad {\cal W}^j \, \in\, \{1,\hdots ,d\}^m \quad , \quad 1\leq j\leq d^m \end{equation}

		\noindent where the words~${\cal W}^j$ are numbered in the lexicographic order. Two control volumes~$C^j_m$ and~$C^l_m$,~$l\neq j$, are said to be \textbf{neighbors}, which we will denote by~$j\underset{m}{\sim}l$, if they have a vertex in common, i.e. if they are linked in the graph~$\mathcal{SS}_m$. A control volume~$C^j_m$ has~$d$ neighbors, or~$d-1$ neighbors if it contains a point of~$V_0$.

	\end{definition}

	\vskip 1cm

	\begin{remark}
		One may check that:

		\begin{equation} \underset{j=1}{\overset{d^m }{\bigcup}}\, C^j_m =   \mathfrak{SS} \, \cdot  \end{equation}

		\noindent Two distinct control volumes intersect at most at one point, which is of measure zero, and:

		$$\forall\,j\,\in\,\{1,\hdots,d^m\}\,:\quad \mu\left(C^j_m\right)=d^{-m}\, \cdot$$

	\end{remark}

	\vskip 1cm
	\noindent We define then, for any integer~$j$ in~$\{1,\hdots,d^m\}$, and any~$t$ in~$[0,T]$, the mean value of~$u$ over the control volume~$C^j_m$:

	\begin{equation}
	u_j(t)=\frac{1}{\mu(C^j_m)}\int_{C^j_m}u(t,X)\, d\mu(X) \, \cdot
	\end{equation}

	\noindent and set:

	\begin{align}
	u^0_j&=u_j(0)=\frac{1}{\mu(C^j_m)}\int_{C^j_m}g(X)\,d\mu(X) \\
	u^n_j&=u_j(t_n)=\frac{1}{\mu(C^j_m)}\int_{C^j_m}u(t_n,X)\, d\mu(X)
	\end{align}

	\noindent The local Gauss-Green formula (see~\cite{StrichartzLivre2006}, (2.3.13)), applied with~$v=1$ on~$C^j_m$, enables one to write:

	\begin{equation} \displaystyle  \int_{C^j_m} \Delta u  \, d\mu = \displaystyle\sum_{X\in\partial C^j_m} \,\partial_n u(X)
  \end{equation}

	\noindent where~$\partial C^j_m=f_{{\cal W}^j}(V_0)$, and where~$\partial_n u(X)$ denotes the local normal derivative of~$u$ at~$X$ relative to the cell~$C^j_m$. Given a natural integer~$n$, one may then integrate the heat equation over $C^j_m\times\left]t_n,t_{n+1}\right[$:

	\begin{equation}
	\label{eq:FVintegrated}
	\displaystyle  \int_{C^j_m} \left(u(t_{n+1},X)-u(t_n,X)\right) d\mu(X) = \displaystyle \int^{t_{n+1}}_{t_n} \sum_{X\in\partial C^j_m} \partial_n u(t,X) \, dt \, \cdot  \end{equation}

	\noindent We recall that $C^j_m=f_{{\cal W}^j}(\mathfrak{SS})$. Its boundary points are the~$d$ points~$X=f_{{\cal W}^j}(P_i)$,~$1\leq i\leq d$. Any such point which does not belong to~$V_0$ is a junction point with exactly one neighbor cell~$C^l_m$:

		\begin{equation}\text{for some }\,
	(k,l)\, \in\, \{1,\hdots ,d \}\times \{1,\hdots ,d^m\} \, :
	\quad  f_{{\cal W}^j}(P_i)=f_{{\cal W}^{l}}(P_k) \, \cdot  \end{equation}

	\noindent One may use the approximation:

	\begin{align}
	\partial_n u(t,X) & \approx \displaystyle \frac{1}{r^m}\sum_{\substack{Y \underset{m}{\sim} X \\ Y\in C^j_m}} \left(u(t,X)-u(t,Y)\right) \\
	&= \displaystyle \frac{1}{r^m}\left((d-1)\,u(t,X)-\sum_{\substack{Y \underset{m}{\sim} X \\ Y\in C^j_m}} u(t,Y)\right)\\
	&= \displaystyle \frac{1}{r^m}\left(d\,u(t,X)-\sum_{Y\in \partial C^j_m} u(t,Y)\right)\\
	&\approx \displaystyle \frac{1}{r^m}\, d \left(u(t,X)-u_j(t)\right)
	\end{align}


	\noindent in conjunction with:

\begin{equation}\displaystyle
	\frac{1}{d}\sum_{Y\in \partial C^j_m} u(t,Y)  \approx \frac{1}{\mu(C^j_m)}\int_{C^j_m}u(t,X)\, d\mu(X)
	= u_j(t) \, \cdot   \end{equation}

	\begin{remark}
		Both approximations are exact when~$u(t,\cdot)$ is a spline function of level~$m$, i.e. belongs to~${\cal S}({\cal H}_0,V_m)$: on each~$m^{th}$-order cell, the local normal derivative of a harmonic function is given by its discrete value at level~$m$, and the mean value of a harmonic function over a cell is the mean value of its values at the~$d$ vertices of the cell (since the standard measure gives the same weight to the~$d$ sub-cells).
	\end{remark}

	\noindent \textbf{Interior boundary points.} At a junction point~$X=f_{{\cal W}^j}(P_i)=f_{{\cal W}^{l}}(P_k)$,~$X\notin V_0$, we then introduce the matching condition:

\begin{equation} \partial_n u(t,f_{{\cal W}^j}(P_i)) =- \partial_n u(t,f_{{\cal W}^{l}}(P_k))  \end{equation}

	\noindent i.e.

	\begin{align}
	\displaystyle \frac{1}{r^m}\, d \left(u(t,X)-u_j(t)\right)
	&= -  \displaystyle \frac{1}{r^m}\, d \left(u(t,X)-u_{l}(t)\right) \, \cdot
	\end{align}

	\noindent This implies:

	\begin{align}
	u(t,X)\approx \frac{u_j(t) + u_{l}(t)}{2} \, \cdot
	\end{align}

	\noindent The normal derivative is then:

	\begin{align}
	\partial_n u(t,X) &\approx \displaystyle \frac{1}{r^m}\,d \left( \frac{u_j(t) + u_{l}(t)}{2}-u_j(t)\right)\\
	&= r^{-m}\, \displaystyle  \frac{d}{2} \left(u_{l}(t) - u_j(t)\right).
	\end{align}

	\noindent \textbf{Points of~$V_0$.} If~$X=f_{{\cal W}^j}(P_i)$ belongs to~$V_0$ (in which case~$X=P_i$, and~${\cal W}^j=(i,\hdots,i)$), the Dirichlet condition~$u(t,X)=0$ yields:

	\begin{align}
	\partial_n u(t,X) &\approx \displaystyle \frac{1}{r^m}\,d \left( 0-u_j(t)\right)=- r^{-m}\, d\, u_j(t)\, \cdot
	\end{align}

	\noindent In other words, the boundary flux is obtained by introducing, beyond each point of~$V_0$, a ghost cell where the mean value is~$-u_j(t)$.\\

	\noindent Back to equation~(\ref{eq:FVintegrated}), with the approximation

	$$\int^{t_{n+1}}_{t_n} \sum_{X\in\partial C^j_m} \partial_n u(t,X) \, dt\approx h\sum_{X\in\partial C^j_m} \partial_n u(t_n,X)$$

	\noindent i.e.

\begin{equation}  u^{n+1}_j\approx u^n_j + \frac{h}{\mu(C^j_m)} \sum_{X\in\partial C^j_m} \partial_n u(t_n,X)   \end{equation}

	\noindent one may now build the finite volume scheme, where we keep, for the sake of simplicity, the notation~$u^n_j$ for the approximate values:

\begin{equation}
\label{eq:FVscheme}
u^{n+1}_j= u^n_j + \frac{h}{\mu(C^j_m)}\, \displaystyle \frac{1}{r^m}\, \frac{d}{2} \left(\sum_{l \underset{m}{\sim}j} \left(u^n_{l} - u^n_j\right) - 2\,\varepsilon_j\, u^n_j\right)   \end{equation}

	\noindent where~$l\underset{m}{\sim}j$ means that the cells $f_{\mathcal{W}^j}(\mathfrak{SS})$ and $f_{\mathcal{W}^l}(\mathfrak{SS})$ are neighbors, and where~$\varepsilon_j=1$ if~$C^j_m$ contains a point of~$V_0$, and~$\varepsilon_j=0$ otherwise. Since~$\mu(C^j_m)=d^{-m}$ and~$r=\displaystyle\frac{d}{d+2}$:

	$$\frac{1}{\mu(C^j_m)}\, \displaystyle \frac{1}{r^m}\, \frac{d}{2}=\frac{d}{2}\,(d+2)^m \, \cdot$$

	\vskip 1cm

	\begin{remark}$\,$\\
		\begin{enumerate}
			\item[\emph{i}.] One may note that we recover, from a totally independent approach, a scheme with the same structure and the same normalization constant~$\displaystyle\frac{d}{2}\,(d+2)^m$ as the finite difference scheme. The two schemes are nevertheless different: the unknowns of the finite volume scheme are the~$d^m$ mean values over the~$m^{th}$-order cells, and the scheme is built on the graph~$\mathcal{SS}_m$, while the~$\mathcal{N}_m-d$ unknowns of the finite difference scheme are the values at the points of~$V_m\setminus V_0$, on the graph~$\mathfrak{SS}_m$.
			\item[\emph{ii}.] We can also define the backward (implicit) scheme:

		\begin{equation}  u^{n}_j= u^{n-1}_j + \frac{h}{\mu(C^j_m)}\, \displaystyle \frac{1}{r^m}\, \frac{d}{2} \left(\sum_{l \underset{m}{\sim}j} \left(u^n_{l} - u^n_j\right) - 2\,\varepsilon_j\, u^n_j\right)  \quad , \quad 1\leq n\leq N \, \cdot \end{equation}
		\end{enumerate}
	\end{remark}

	\vskip 1cm

	\noindent Given~$m\,\in \,\N^\star$, we number the control volumes~$C^j_m$,~$1\leq j\leq d^m$, according to the lexicographic order of the words~${\cal W}^j$, and consider:

	\begin{equation} \forall\, n \, \in\, \left \lbrace 0, \hdots, N \right \rbrace\, :\quad
	U(n) =\left(
	\begin{matrix}
	u^n_1\\
	\vdots\\
	u^n_{d^m}\\
	\end{matrix}
	\right) \end{equation}

	\noindent (one has to bear in mind that there are $d^m$~$m^{th}$-order cells).\\
	\noindent One has:

	$$\forall\, n \, \in\, \left \lbrace 0, \hdots, N-1 \right \rbrace\, :\quad U(n+1) = A \, U(n)$$

	\noindent where:

\begin{equation}
	A=I_{d^m}-h\,\tilde{\Delta}_{m} \quad , \quad \tilde{\Delta}_{m}=\frac{d}{2}\,(d+2)^m\,\tilde{L}_m
	\end{equation}

	\noindent and where~$I_{d^m}$ denotes the~$ {d^m}\times  {d^m}$ identity matrix, and~$\tilde{L}_m$ the~$ {d^m}\times  {d^m}$ matrix of the graph Laplacian of~$\mathcal{SS}_m$, with Dirichlet boundary terms: its diagonal coefficient~$(j,j)$ is equal to~$d$ if~$C^j_m$ does not contain any point of~$V_0$, and to~$(d-1)+2=d+1$ otherwise, and, for~$j\neq l$, its coefficient~$(j,l)$ is equal to~$-1$ if~$j\underset{m}{\sim}l$, and to~$0$ otherwise. The matrix~$\tilde{L}_m$ is symmetric, positive definite.

	\vskip 1cm
	
\subsubsection{Consistency, stability and convergence}

	\paragraph{Theoretical study of the error}
	$\,$\\

	\noindent In the sequel, we make the following assumptions on the solution~$u$ of problem~(\ref{eq:HeatFVM}):

	\begin{itemize}
		\item[(H1)] $\displaystyle\frac{\partial u}{\partial t}$ and~$\displaystyle\frac{\partial^2 u}{\partial t^2}$ are continuous and bounded on~$[0,T]\times\mathfrak{SS}$;
		\item[(H2)] for any~$t\in[0,T]$,~$u(t,\cdot)\,\in\,\text{dom}\,\Delta$, and~$\Delta u(t,\cdot)$ is H\"{o}lder continuous of exponent~$\alpha$, uniformly with respect to~$t$: there exists a constant~$C_\Delta>0$ such that
		$$\forall\, (X,Y)\,\in \mathfrak{SS}^2,\,\forall\,t\in[0,T] : \quad |\Delta u(t, X)-\Delta u(t, Y)| \leq C_\Delta  \, |X-Y|^\alpha$$
		where~$\alpha=\displaystyle\frac{\ln\left(\frac{d+2}{d}\right)}{2\,\ln 2}$ is the H\"{o}lder exponent of the functions of finite energy (this assumption is satisfied as soon as~$\Delta u(t,\cdot)$ is of finite energy, uniformly with respect to~$t$).
	\end{itemize}

	\noindent Since~$u(t,P)=0$ for any~$t$ and any~$P\,\in\,V_0$, assumption~(H1) yields~$\displaystyle\frac{\partial u}{\partial t}(t,P)=0$, and thus, by continuity:

	\begin{equation}
	\label{eq:CompatibilityFVM}
	\forall\,t\,\in\,[0,T],\,\forall\,P\,\in\,V_0\,:\quad \Delta u(t,P)=0 \, \cdot
	\end{equation}

	\noindent We denote by~$\delta$ the diameter of~$\mathfrak{SS}$, so that the diameter of any~$m^{th}$-order cell is~$2^{-m}\,\delta$.\\

	\noindent \textbf{Time discretization.} Under assumption~(H1), for any~$n\,\in\,\{0,\hdots,N-1\}$ and any~$j\,\in\,\{1,\hdots,d^m\}$:

	\begin{equation}
	\frac{u_j(t_{n+1})-u_j(t_n)}{h}=\frac{1}{\mu(C^j_m)}\int_{C^j_m}\frac{\partial u}{\partial t}(t_n,X)\,d\mu(X)+{\mathcal O}(h)=\frac{1}{\mu(C^j_m)}\int_{C^j_m}\Delta u(t_n,X)\,d\mu(X)+{\mathcal O}(h) \, \cdot
	\end{equation}

	\noindent Let us note that the local error over one time step is of order~${\mathcal O}(h^2)$, which, after division by~$h$, yields a consistency error of order~${\mathcal O}(h)$, as for the explicit finite difference scheme.\\

	\noindent \textbf{Space discretization.} We follow the method introduced by R.~S.~Strichartz in the case of the Gasket (see~\cite{StrichartzAverages2001}, Theorem~4.1, where the uniform convergence of the renormalized cell graph Laplacian~$\frac{3}{2}\,5^m\left(\sum_{w'\sim w}A_{w'}(u)-3\,A_w(u)\right)$ towards~$\Delta u$ is proved, for non-boundary cells). We extend it to Sierpi\'nski simplices, include the boundary cells thanks to the Dirichlet terms, and obtain a rate of convergence under assumption~(H2). We need two lemmas.

	\begin{lemma}\textbf{(Discrete identity)}\\
		\label{lem:DiscreteIdentity}
		Let~$m\,\in\,\N^\star$, and~$v$ a real-valued function defined on~$V_m$, which vanishes on~$V_0$. For any~$j\,\in\,\{1,\hdots,d^m\}$, let us set:
		$$\bar{v}_j=\frac{1}{d}\sum_{X\in\partial C^j_m}v(X)\, \cdot$$
		Then, for any~$j\,\in\,\{1,\hdots,d^m\}$:
		\begin{equation}
		\sum_{l \underset{m}{\sim}j}\left(\bar{v}_l-\bar{v}_j\right)-2\,\varepsilon_j\,\bar{v}_j=\frac{1}{d}\sum_{X\in\partial C^j_m\setminus V_0}\Delta_m v(X) \, \cdot
		\end{equation}
	\end{lemma}

	\begin{proof}
		Let~$X\,\in\,\partial C^j_m\setminus V_0$, and~$C^l_m$ the neighbor cell of~$C^j_m$ at~$X$. Then:
		$$d\left(\bar{v}_l-\bar{v}_j\right)=\sum_{\substack{Y\in \partial C^l_m\\ Y\neq X}}\left(v(Y)-v(X)\right)-\sum_{\substack{Y\in \partial C^j_m\\ Y\neq X}}\left(v(Y)-v(X)\right)
		=\Delta_m v(X)-2\sum_{\substack{Y\in \partial C^j_m\\ Y\neq X}}\left(v(Y)-v(X)\right) \, \cdot$$
		If~$C^j_m$ contains a point~$P$ of~$V_0$, one has, since~$v(P)=0$:
		$$-2\,d\,\bar{v}_j=-2\sum_{\substack{Y\in \partial C^j_m\\ Y\neq P}}\left(v(Y)-v(P)\right) \, \cdot$$
		By summation over all the points~$X$ of~$\partial C^j_m$, the terms~$\displaystyle\sum_{X\in\partial C^j_m}\,\sum_{\substack{Y\in \partial C^j_m\\ Y\neq X}}\left(v(Y)-v(X)\right)$ vanish by antisymmetry, which yields the result.
	\end{proof}

	\vskip 0.5cm

	\begin{lemma}\textbf{(Mean value over a cell)}\\
		\label{lem:MeanValue}
		Let us denote by~$\zeta$ the unique function of~$\text{dom}\,\Delta$ such that~$-\Delta\zeta=1$ on~$\mathfrak{SS}$, and~$\zeta_{\vert V_0}=0$, and set:
		$$\kappa=\int_{\mathfrak{SS}}\zeta\,d\mu>0 \, \cdot$$
		Let~$v\,\in\,\text{dom}\,\Delta$, such that~$\Delta v$ is H\"{o}lder continuous of exponent~$\alpha$, with constant~$C_\Delta$. Then, for any~$m^{th}$-order cell~$C^j_m$, and any point~$c$ of~$C^j_m$:
		\begin{equation}
		\frac{1}{\mu(C^j_m)}\int_{C^j_m}v\,d\mu=\frac{1}{d}\sum_{X\in\partial C^j_m}v(X)-\kappa\,(d+2)^{-m}\,\Delta v(c)+{\mathcal O}\left((d+2)^{-m}\,2^{-m\alpha}\right)
		\end{equation}
		where the constant in~${\mathcal O}$ only depends on~$\kappa$,~$C_\Delta$ and~$\delta$.
	\end{lemma}

	\begin{proof}
		Let us write~$C^j_m=f_{\cal W}(\mathfrak{SS})$, with~$|{\cal W}|=m$, and set~$w=v\circ f_{\cal W}$. The self-similarity of the measure~$\mu$ yields:
		$$\frac{1}{\mu(C^j_m)}\int_{C^j_m}v\,d\mu=\int_{\mathfrak{SS}}w\,d\mu \, \cdot$$
		The scaling property of the Laplacian (see~\cite{StrichartzLivre2006}) yields~$\Delta w=r^m\,d^{-m}\,(\Delta v)\circ f_{\cal W}=(d+2)^{-m}\,(\Delta v)\circ f_{\cal W}$.\\
		Let~$h_w$ denote the harmonic function which takes the same values as~$w$ on~$V_0$. Since the standard measure is invariant under the symmetries of the simplex, the integral of a harmonic function is the mean value of its values on~$V_0$:
		$$\int_{\mathfrak{SS}}h_w\,d\mu=\frac{1}{d}\sum_{i=1}^d w(P_i)=\frac{1}{d}\sum_{X\in\partial C^j_m}v(X) \, \cdot$$
		The function~$w-h_w$ vanishes on~$V_0$, and~$\Delta(w-h_w)=\Delta w$. Thanks to the Green function~$G$ of the Dirichlet problem (see~\cite{StrichartzLivre2006}, section~2.6):
		$$w-h_w=-\int_{\mathfrak{SS}}G(\cdot,Y)\,\Delta w(Y)\,d\mu(Y)\, , \quad \text{and}\quad \zeta=\int_{\mathfrak{SS}}G(\cdot,Y)\,d\mu(Y) \, \cdot$$
		Integrating, and using the symmetry of~$G$:
		$$\int_{\mathfrak{SS}}(w-h_w)\,d\mu=-(d+2)^{-m}\int_{\mathfrak{SS}}\zeta(Y)\,(\Delta v)(f_{\cal W}(Y))\,d\mu(Y) \, \cdot$$
		Since~$\zeta\geq 0$, and since~$|(\Delta v)(f_{\cal W}(Y))-\Delta v(c)|\leq C_\Delta\,(2^{-m}\,\delta)^\alpha$ for any~$Y\,\in\,\mathfrak{SS}$, one gets the result.
	\end{proof}

	\vskip 0.5cm

	\begin{proposition}\textbf{(Space consistency)}\\
		\label{prop:SpaceConsistency}
		Under assumption~(H2), for any~$t\,\in\,[0,T]$ and any~$j\,\in\,\{1,\hdots,d^m\}$:
		\begin{equation}
		\frac{d}{2}\,(d+2)^m\left(\sum_{l \underset{m}{\sim}j}\left(u_l(t)-u_j(t)\right)-2\,\varepsilon_j\,u_j(t)\right)=\frac{1}{\mu(C^j_m)}\int_{C^j_m}\Delta u(t,X)\,d\mu(X)+{\mathcal O}\left(2^{-m\alpha}\right)
		\end{equation}
		uniformly with respect to~$t$ and~$j$.
	\end{proposition}

	\begin{proof}
		Let us fix~$t$, and apply lemma~\ref{lem:MeanValue} to~$v=u(t,\cdot)$, on the cell~$C^j_m$ and on its neighbors~$C^l_m$:
		$$u_k(t)=\bar{u}_k-\kappa\,(d+2)^{-m}\,\Delta u(t,c_k)+{\mathcal O}\left((d+2)^{-m}\,2^{-m\alpha}\right)\quad , \quad \bar{u}_k=\frac{1}{d}\sum_{X\in\partial C^k_m}u(t,X)$$
		where, for any~$k$,~$c_k$ is a point of~$C^k_m$. For any neighbor~$C^l_m$ of~$C^j_m$, we choose for~$c_l$ the junction point of~$C^j_m$ and~$C^l_m$; we choose~$c_j=P$ if~$C^j_m$ contains a point~$P$ of~$V_0$, and any point of~$C^j_m$ otherwise. All the points~$c_l$ and~$c_j$ thus belong to~$C^j_m$.\\

		\noindent \emph{i}. Contribution of the vertex averages. Lemma~\ref{lem:DiscreteIdentity}, applied to~$v=u(t,\cdot)_{\vert V_m}$, yields:
		$$\frac{d}{2}\,(d+2)^m\left(\sum_{l \underset{m}{\sim}j}\left(\bar{u}_l-\bar{u}_j\right)-2\,\varepsilon_j\,\bar{u}_j\right)=\frac{1}{d}\sum_{X\in\partial C^j_m\setminus V_0}\frac{d}{2}\,(d+2)^m\,\Delta_m u(t,X) \, \cdot$$
		Under assumption~(H2), the pointwise estimate of the finite difference method (see~\cite{StrichartzLivre2006}, exercise~2.4.1) yields~$\displaystyle\frac{d}{2}\,(d+2)^m\,\Delta_m u(t,X)=\Delta u(t,X)+{\mathcal O}(2^{-m\alpha})$. Thanks to~(\ref{eq:CompatibilityFVM}), the points of~$V_0$ can be added to the sum, which yields:
		$$\frac{1}{d}\sum_{X\in\partial C^j_m}\Delta u(t,X)+{\mathcal O}(2^{-m\alpha})=\frac{1}{\mu(C^j_m)}\int_{C^j_m}\Delta u(t,X)\,d\mu(X)+{\mathcal O}(2^{-m\alpha})$$
		since the oscillation of~$\Delta u(t,\cdot)$ over~$C^j_m$ is bounded by~$C_\Delta\,(2^{-m}\,\delta)^\alpha$.\\

		\noindent \emph{ii}. Contribution of the correction terms. They write:
		$$-\frac{d}{2}\,\kappa\left(\sum_{l \underset{m}{\sim}j}\left(\Delta u(t,c_l)-\Delta u(t,c_j)\right)-2\,\varepsilon_j\,\Delta u(t,c_j)\right)+{\mathcal O}(2^{-m\alpha})={\mathcal O}(2^{-m\alpha})$$
		since the points~$c_l$ and~$c_j$ belong to~$C^j_m$, and are thus at distance at most~$2^{-m}\,\delta$, and since, if~$\varepsilon_j=1$,~$\Delta u(t,c_j)=\Delta u(t,P)=0$, thanks to~(\ref{eq:CompatibilityFVM}).
	\end{proof}

	\vskip 1cm

	\paragraph{Consistency}

	\begin{definition}
		The consistency error of the scheme~(\ref{eq:FVscheme}) is defined, for~$0\leq n \leq N-1$ and~$1\leq j \leq d^m$, by:
		$$\varepsilon^m_{n,j}=\frac{u_j(t_{n+1})-u_j(t_n)}{h}-\frac{d}{2}\,(d+2)^m\left(\sum_{l \underset{m}{\sim}j}\left(u_l(t_n)-u_j(t_n)\right)-2\,\varepsilon_j\,u_j(t_n)\right)$$
		where the~$u_j(t)$ denote the mean values of the exact solution. The scheme is said to be \textbf{consistent} if the consistency error tends towards zero when $h\rightarrow 0$ and~\mbox{$m \rightarrow+ \infty$}, for a given norm.
	\end{definition}

	\noindent For $0\leq n \leq N-1,\, 1\leq j \leq d^m$, the time discretization estimate and proposition~\ref{prop:SpaceConsistency} yield, under assumptions~(H1) and~(H2):

	\begin{equation}
	\varepsilon^m_{n,j} = \mathcal{O}(h) + \mathcal{O}\left(2^{-m\alpha}\right)
	\end{equation}

	\noindent uniformly with respect to~$n$ and~$j$. One has:

\begin{equation}  \displaystyle \lim_{h \to 0^+, \, m \to+\infty}\,\max_{n,j}\, |\varepsilon^m_{n,j}|=0	\end{equation}

	\noindent which yields the consistency of the scheme.

	\vskip 1cm

	\paragraph{Stability}

	\begin{notation}\textbf{Spectral radius}\\
		Given a square matrix~$A$, we will denote by~$\lambda_{\max}(A)$ its spectral radius.
	\end{notation}

	\vskip 1cm

	\begin{notation}\textbf{Spectral norm}\\
		Given a square matrix~$A$, we will denote by~$\rho(A) $ its spectral norm, obtained through:

		\begin{equation} \rho(A) = \left ( \lambda_{\max} \, \left( A^T \, A\right ) \right )^{\frac{1}{2}} 	\end{equation}

		\noindent For a symmetric matrix~$A$,~$\rho(A)$ is the largest absolute value of the eigenvalues of~$A$.
	\end{notation}

	\vskip 1cm

	\begin{proposition}\textbf{(Spectral decimation on the cell graphs)}\\
		\label{prop:DecimationFVM}
		
		\noindent We introduce the real valued function~$\Phi$, defined on~$\R$, by:
		
		$$ \forall\, t \,\in\, \R :\quad  \Phi(t)= (d+2-t)\,t\,  \cdot$$
		\noindent For any integer~$m\geq 2$, and any eigenvalue~$\lambda_{m}$ of the matrix~$\tilde{L}_m$ which is different from~$d$ and from~$d+2$, $\Phi(\lambda_m)$ is an eigenvalue of~$\tilde{L}_{m-1}$:
		
		\begin{equation}  \lambda_{m-1}=\Phi \left (\lambda_{m}\right ) \, \cdot 	\end{equation}
		
		\noindent The values~$d$ and~$d+2$ are created at each level. Moreover, the eigenvalues of~$\tilde{L}_1$ are~$2$, with multiplicity~$1$, and~$d+2$, with multiplicity~$d-1$.
		
	\end{proposition}
	
	\vskip 1cm
	
	\begin{proof}
		\noindent In the case of the Gasket, this result is due to R.~J.~Ravier and R.~S.~Strichartz (see~\cite{RavierStrichartz2015}, Theorem~3.2): an eigenfunction~$u$ of the cell graph at level~$m+1$, with eigenvalue~$\lambda_{m+1}$, yields, by averaging over the~$d$ sub-cells of each cell of level~$m$,
		$$u(X_k)=\displaystyle \frac{1}{d}\sum_{i=1}^{d} u(a^k_i)\, ,$$
		an eigenfunction at level~$m$, with eigenvalue~$\lambda_m=\lambda_{m+1}\,(d+2-\lambda_{m+1})$. The general case follows from the same computation. Let us consider the sequence of graphs~$(\mathcal{SS}_m)_{m\geq 1}$, and a vertex~$a^k_i$ of~$\mathcal{SS}_{m+1}$, i.e. a sub-cell of the~$k^{th}$ cell~$X_k$ of level~$m$, of neighbors~$a^k_j$,~$j\neq i$, and~$a^l_h$, which belongs to another cell~$X_l$ of level~$m$ (see figures~\ref{figSpecDec0} and~\ref{figSpecDec1}). The eigenvalue relation writes:
		
		\begin{equation} \left \lbrace d-\lambda_{m+1}\right \rbrace \, u(a^k_i)=\displaystyle\sum_{j\neq i} u(a^k_j)+u(a^l_h) \, \cdot	\end{equation}
		\begin{figure}[!htb]
			
			\minipage{0.29\textwidth}
			\includegraphics[width=\linewidth]{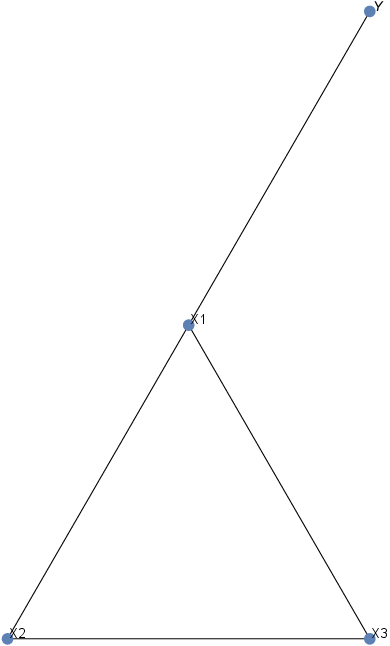}
			\caption{$\mathcal{SS}_m$ for the Sierpi\'{n}ski triangle.}
			\label{figSpecDec0}
			\endminipage\hfill
			\minipage{0.4\textwidth}
			\includegraphics[width=\linewidth]{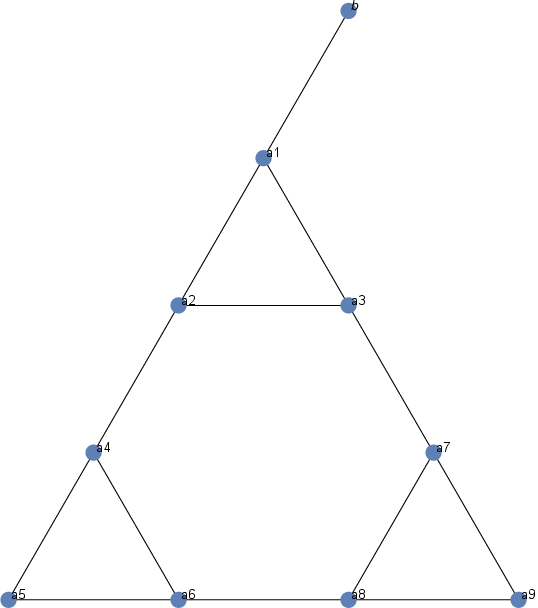}
			\caption{$\mathcal{SS}_{m+1}$ for the Sierpi\'{n}ski triangle.}
			\label{figSpecDec1}
			\endminipage\hfill
		\end{figure}
		
		\noindent Using the mean property:
		
	\begin{equation} u(X_k)=\displaystyle \frac{1}{d}\sum_{i=1}^{d} u(a^k_i) 	\end{equation}
		
		\noindent one obtains, by adding~$u(a^k_i)$ to both sides of the eigenvalue relation:
		
	\begin{equation} (d+1-\lambda_{m+1})\,u(a^k_i)=d \, u(X_k)+u(a^l_h)\quad , \quad (d+1-\lambda_{m+1})\,u(a^l_h)=d \, u(X_l)+u(a^k_i) 	\end{equation}
		
		\noindent which leads, for~$\lambda_{m+1}\neq d$ and~$\lambda_{m+1}\neq d+2$, to:
		
	\begin{equation} u(a^k_i)=d\, \displaystyle \frac{((d+1)-\lambda_{m+1})\,u(X_k)+u(X_l)}{(d+2-\lambda_{m+1})(d-\lambda_{m+1})} \, \cdot 	\end{equation}
		
		\noindent The sub-cells which contain a point of~$V_0$ are treated in the same way, the Dirichlet term playing the role of a neighbor where the mean value is~$-u$. By summation over the~$d$ sub-cells~$a^k_i$ of~$X_k$, one obtains the eigenvalue relation at~$X_k$, with:
		
			\begin{equation} \lambda_m=\lambda_{m+1}\, \left(d+2-\lambda_{m+1}\right) \, \cdot 	\end{equation}
		
		\noindent Finally,~$\tilde{L}_1=(d+2)\,I_d-J_d$, where~$J_d$ denotes the~$d\times d$ matrix the coefficients of which are all equal to~$1$, the eigenvalues of which are~$d$, with multiplicity~$1$, and~$0$, with multiplicity~$d-1$.
	\end{proof}

	\vskip 1cm
	
	\noindent Thus:
	
	\begin{equation}  \forall\, m\,\geq\,2\, : \quad 
	\left \lbrace 
	\begin{array}{ccc} \lambda^{-}_{m}&=&\displaystyle\frac{(d+2) -
		\left ( (d+2)^2-4\, \lambda_{m-1}\right)^{\frac{1}{2}}}{2}\\
	\lambda^{+}_{m}&=&\displaystyle\frac{(d+2) +
		\left ( (d+2)^2-4\, \lambda_{m-1}\right)^{\frac{1}{2}}}{2}\end{array} \right. 
	\, \cdot 	\end{equation}
	
	\noindent It is then natural to consider the maps~$\psi^-$ and~$\psi^+$ defined as follows:

		\begin{equation}   \forall\, x\,\in\, \left]-\infty,\displaystyle\frac{(d+2)^2}{4}\right] \, : \quad 
	\left \lbrace 
	\begin{array}{ccc}\psi^-(x)& =&\displaystyle\frac{(d+2)- \left ( (d+2)^2-4\, x\right)^{\frac{1}{2}}}{2}\\
	\psi^+(x)&=&\displaystyle\frac{(d+2) +  \left ( (d+2)^2-4\, x\right)^{\frac{1}{2}}}{2}
	\end{array}\right. \,\cdot
	\end{equation}

	\noindent One has:
	
	\begin{equation}  \left \lbrace\begin{array}{ccc}\psi^-(0)&=&0 \\  \psi^-\left(\displaystyle\frac{(d+2)^2}{4}\right)&=&\displaystyle\frac{d+2}{2}\\
	\end{array}
	\right. \quad , \quad 
	\left \lbrace\begin{array}{ccc}\psi^+(0)&=&d+2 \\  
	\psi^+\left(\displaystyle\frac{(d+2)^2}{4}\right)&=&\displaystyle\frac{d+2}{2}
	\end{array}
	\right. 	\end{equation}
	
	\noindent In particular,~$\psi^-$ and~$\psi^+$ take their values in~$[0,d+2]$. The maps~$\psi^-$ and~$\psi^+$ are respectively non decreasing, and non increasing, with respective fixed points:

		\begin{equation}  \left \lbrace\begin{array}{ccc}
	x^{-,\star}&=&0\\
	x^{+,\star}&=&(d+2)-1\end{array}\right. 	\end{equation}
	
	\noindent Another interesting feature of those two maps is their contractive property in a neighbourhood of their fixed points, due to:

	\begin{equation}  \left \lbrace\begin{array}{ccccccc}
	\left |\displaystyle \frac{d}{dx}\psi^-(0)\right | &=&\displaystyle\frac{1}{\sqrt{(d+2)^2}}
	&=&\displaystyle\frac{1}{d+2}&<&1 \\
	\left |\displaystyle \frac{d}{d\, x}\psi^+\left ((d+2)-1\right )\right |&=&\displaystyle\frac{1}{\sqrt{(d+2)^2-4\, (d+2)+4}} 
	&=&\displaystyle\frac{1}{d}&<&1
	\end{array}\right. 	\end{equation}
	
	\noindent The whole spectrum of~$\tilde{L}_m$, for~$m\geq 2$, is thus generated by the maps~$\psi^-$ and~$\psi^+$, applied to the eigenvalues of the previous level, together with the values~$d$ and~$d+2$ created at level~$m$. Since~$\psi^\pm$ take their values in~$[0,d+2]$, this yields the bound~$\lambda_m\,\in\,]0,d+2]$, which is also obtained, by a different argument, in proposition~\ref{prop:SpectralInclusion} below. We first give an elementary bound, which is sufficient for the stability of the scheme.
	
	\vskip 1cm
	
	\begin{proposition}\textbf{(Spectral bound)}\\
		\label{prop:SpectralBoundFVM}
		Any eigenvalue~$\lambda$ of the matrix~$\tilde{L}_m$ satisfies:

		\begin{equation}
		0<\lambda\leq 2\,d \, \cdot
		\end{equation}

		\noindent Moreover, the largest eigenvalue~$\lambda_{\max}(\tilde{L}_m)$ satisfies~$\lambda_{\max}(\tilde{L}_m)\geq d+1$.
	\end{proposition}

	\vskip 1cm

	\begin{proof}
		The matrix~$\tilde{L}_m$ is symmetric. Its diagonal coefficient~$(j,j)$ is equal to~$d$ if~$C^j_m$ does not contain any point of~$V_0$, and to~$d+1$ otherwise, while the sum of the absolute values of its off-diagonal coefficients on the~$j^{th}$ row is equal to the number of neighbors of~$C^j_m$, i.e. to~$d$, or to~$d-1$. Gershgorin's theorem yields:~$0\leq\lambda\leq 2\,d$. Since~$\tilde{L}_m$ is diagonally dominant, irreducible, and strictly diagonally dominant on the rows related to the cells which contain a point of~$V_0$, it is positive definite, which yields~$\lambda>0$.\\
		Finally, if~$C^j_m$ contains a point of~$V_0$, and if~$e_j$ denotes the~$j^{th}$ vector of the canonical basis, one has~$\lambda_{\max}(\tilde{L}_m)\geq e_j^T\,\tilde{L}_m\,e_j=d+1$.
	\end{proof}

	\vskip 1cm

	\noindent The bound~$\lambda\leq 2\,d$ can be improved thanks to the following result, which relates the spectrum of~$\tilde{L}_m$ to the one of the finite difference scheme.
	
	\vskip 1cm
	
	\begin{proposition}\textbf{(Comparison with the finite difference spectrum)}\\
		\label{prop:SpectralInclusion}
		Let~$m\,\in\,\N^\star$, and~$d\geq 3$. Let us denote by~$L_m$ the~$(\mathcal{N}_m-d)\times(\mathcal{N}_m-d)$ matrix of the Dirichlet graph Laplacian~$-\Delta_m$ on~$V_m\setminus V_0$ (i.e. the matrix of the finite difference scheme), and by~$\mathbb{A}_m$ the~$d^m\times(\mathcal{N}_m-d)$ averaging matrix, which, to any function~$v$ defined on~$V_m$ and vanishing on~$V_0$, associates the vector of its vertex averages~$\displaystyle\bar{v}_j=\frac{1}{d}\sum_{X\in\partial C^j_m}v(X)$,~$1\leq j\leq d^m$. Then:
		\begin{enumerate}
			\item[\emph{i}.] $\mathbb{A}_m\,L_m=\tilde{L}_m\,\mathbb{A}_m$;
			\item[\emph{ii}.] $\mathbb{A}_m$ is onto;
			\item[\emph{iii}.] the spectrum of~$\tilde{L}_m$ is contained in the Dirichlet spectrum of~$L_m$;
			\item[\emph{iv}.] $2\,d$ is not an eigenvalue of~$\tilde{L}_m$;
			\item[\emph{v}.] any eigenvalue~$\lambda$ of~$\tilde{L}_m$ satisfies:~$0<\lambda\leq d+2$.
		\end{enumerate}
	\end{proposition}
	
	\begin{proof}
		\emph{i}. This is a reformulation of lemma~\ref{lem:DiscreteIdentity}, since~$L_m\,v=-\Delta_m v$ on~$V_m\setminus V_0$.\\
		
		\noindent \emph{ii}. Let~$w\,\in\,\R^{d^m}$ such that~$\mathbb{A}_m^T\,w=0$. For any~$X\,\in\,V_m\setminus V_0$, which belongs to exactly two cells~$C^j_m$ and~$C^l_m$, one has~$\left(\mathbb{A}_m^T\,w\right)(X)=\displaystyle\frac{1}{d}\left(w_j+w_l\right)=0$. Thus~$w_l=-w_j$ for any pair of neighbor cells. Since~$d\geq 3$, the~$d$ sub-cells of any~$(m-1)^{th}$-order cell are pairwise neighbors, and contain three cells~$C^{j_1}_m$,~$C^{j_2}_m$,~$C^{j_3}_m$ which are pairwise neighbors, which yields~$w_{j_1}=-w_{j_2}=w_{j_3}=-w_{j_1}$, i.e.~$w_{j_1}=0$. Since the graph~$\mathcal{SS}_m$ is connected, one gets~$w=0$:~$\mathbb{A}_m^T$ is one-to-one, and~$\mathbb{A}_m$ is onto.\\
		
		\noindent \emph{iii}. The matrices~$L_m$ and~$\tilde{L}_m$ are symmetric. By transposition of~\emph{i},~$L_m\,\mathbb{A}_m^T=\mathbb{A}_m^T\,\tilde{L}_m$. If~$\tilde{L}_m\,w=\lambda\,w$, with~$w\neq 0$, then~$L_m\left(\mathbb{A}_m^T\,w\right)=\lambda\,\mathbb{A}_m^T\,w$, and~$\mathbb{A}_m^T\,w\neq 0$ thanks to~\emph{ii}: $\lambda$ is an eigenvalue of~$L_m$.\\
		
		\noindent \emph{iv}. Assume that~$\tilde{L}_m\,w=2\,d\,w$, with~$w\neq 0$, and choose~$j$ such that~$|w_j|=\displaystyle\max_k |w_k|$. The~$j^{th}$ row of this relation writes~$\displaystyle\sum_{l \underset{m}{\sim}j}w_l=-\left(2\,d-\tilde{L}_{m,jj}\right)w_j$, where~$2\,d-\tilde{L}_{m,jj}$ is precisely the number of neighbors of~$C^j_m$. This is only possible if~$w_l=-w_j$ for all the neighbors~$C^l_m$ of~$C^j_m$, which then also realize the maximum of~$|w_k|$. By connectedness,~$w_l=-w_j$ for any pair of neighbor cells, which is impossible, as shown in~\emph{ii}.\\
		
		\noindent \emph{v}. According to~\cite{Fukushima1992} (Proposition~3.2), the Dirichlet eigenvalues of~$L_m$ which are different from~$2\,d$ are either equal to~$2$ or to~$d+2$, or are images of eigenvalues of the previous level by one of the maps~$\phi^-$ or~$\phi^+$, which take their values in~$[0,d+2]$. They thus belong to~$]0,d+2]$, and the result follows from~\emph{iii} and~\emph{iv}.
	\end{proof}
	
	\vskip 0.5cm
	
	\begin{remark}
		In the case of the Gasket, the fact that average values of eigenfunctions on the graphs~$\mathfrak{SS}_m$ yield eigenfunctions on the cell graphs has been noticed by R.~J.~Ravier (see~\cite{Ravier}, Lemma~4.4). Numerically, one observes that~$\lambda_{\max}(\tilde{L}_m)=d+2$ for any~$m\geq 1$, and that the smallest eigenvalues of~$\tilde{L}_m$ and~$L_m$ coincide.
	\end{remark}

	\vskip 1cm

	\begin{notation}\textbf{Matrix norm}\\
		In the sequel, we will work with the matrix norm~$\|\cdot \|$ defined, for any square matrix~$A$, by:\\

		$$\| A\| = \parallel A \parallel_2=\rho(A)$$

	\end{notation}

	\vskip 1cm

	\begin{definition}{\ }\\

		Let us set~$\delta x=2^{-m}$. We will say that:

		\begin{enumerate}
			\item[\emph{i}.]  The scheme is~\emph{unconditionally stable} if there exists a real
			constant~\mbox{$C \geq 1$}, independent of~$h$ and~$m$, such that:

			\begin{equation} \forall \, n\, \in \, \left \lbrace 1,\hdots,N\right\rbrace\, :
			\quad \rho(A^{n})\leq C  \, \cdot	\end{equation}

			\item[\emph{ii}.] The scheme is \emph{conditionally stable} if there exist three strictly positive constants~$\beta $, $C_1 $ and $C_2\geq 1$, independent of~$h$ and~$m$, such that:

			\begin{equation}
			\forall \, n\, \in \, \left \lbrace 1,\hdots,N\right\rbrace\, :
			\quad h \leq C_1\,(\delta x)^{\beta}\Longrightarrow \rho(A^{n})\leq C_2  \, \cdot	\end{equation}

		\end{enumerate}
	\end{definition}

	\vskip 1cm

	\begin{proposition}{\ }\\

		\noindent If $\text{Sp}\, (A)=\left \lbrace \gamma_1\, \hdots, \gamma_{d^m}\right \rbrace $ is the spectrum of the matrix~$A$, one has:

		\begin{equation}
		\label{eq:CFLFVM}
		h\, (d+2)^{m}\leq \displaystyle\frac{2}{d^2} \Longrightarrow \forall\,j \,\in\, \left \lbrace  1,\hdots,d^m\right \rbrace \,: \quad  |\gamma_j|\leq 1  \, \cdot	\end{equation}

		\noindent Since~$(d+2)^{-m}=(\delta x)^{D_W}$, where~$D_W=\displaystyle\frac{\ln(d+2)}{\ln 2}$ is the walk dimension of~$\mathfrak{SS}$, the explicit finite volume scheme is conditionally stable, with~$\beta=D_W$,~$C_1=\displaystyle\frac{2}{d^2}$ and~$C_2=1$ (Courant-Friedrichs-Lewy condition). Conversely, if~$h\,(d+2)^{m}>\displaystyle\frac{4}{d\,(d+1)}$, the scheme is unstable: the exponent~$\beta=D_W$ is optimal.
	\end{proposition}

	\vskip 1cm

	\begin{proof}

		\noindent Since:

		\begin{equation}\forall\, n\,\in\,  \left \lbrace 0, \hdots, N-1\right \rbrace \, :
		\quad U(n+1) = A \, U(n) 	\end{equation}

		\noindent where:

	\begin{equation}
		A=I_{d^m}-h\,\tilde{\Delta}_{m}=I_{d^m}-h\,\frac{d}{2}\,(d+2)^m\,\tilde{L}_m
		\end{equation}

		\noindent one has:

		\begin{equation}\forall\, n\,\in\,  \left \lbrace 1, \hdots, N\right \rbrace \, :\quad U(n)=A^{n}\, U(0)  	\end{equation}

		\noindent The matrix~$A$ is symmetric, and its eigenvalues are given by:

	\begin{equation}\forall\, i\,\in\,  \left \lbrace 1, \hdots, d^m\right \rbrace \, :\quad \gamma_i = 1 - h\,  \displaystyle  \frac{d}{2}(d+2)^{m} \, \lambda_i 	\end{equation}

		\noindent where the~$\lambda_i$ denote the eigenvalues of~$\tilde{L}_m$. Thanks to proposition~\ref{prop:SpectralBoundFVM}:

		\begin{equation}\forall\, i\,\in\,  \left \lbrace 1, \hdots, d^m\right \rbrace \, :\quad 1 -   \displaystyle \frac{h\, d}{2}(d+2)^{m}\times  2\, d  \leq\gamma_i< 1 \, \cdot 	\end{equation}

		\noindent Consequently, if the condition:

		\begin{equation}h\, (d+2)^{m}\leq \displaystyle\frac{2}{d^2} 	\end{equation}

		\noindent is satisfied, one obtains~$|\gamma_i|\leq 1$ for any~$i$, and, since~$A$ is symmetric,~$\rho(A^n)=\rho(A)^n\leq 1$.\\

		\noindent Conversely, since~$\lambda_{\max}(\tilde{L}_m)\geq d+1$, the smallest eigenvalue of~$A$ satisfies~$\gamma_{\min}\leq 1-h\,\displaystyle\frac{d}{2}\,(d+2)^m\,(d+1)$, which is strictly less than~$-1$ as soon as~$h\,(d+2)^{m}>\displaystyle\frac{4}{d\,(d+1)}$. In this case,~$\rho(A^n)=\rho(A)^n\to+\infty$.

	\end{proof}
	\vskip 1cm

	\begin{corollary}\textbf{(Sharp stability condition)}\\
		Thanks to proposition~\ref{prop:SpectralInclusion}, the explicit finite volume scheme is stable as soon as:
		$$h\,(d+2)^{m+1}\leq\frac{4}{d}\, ,$$
		i.e.~$h\leq\displaystyle\frac{4}{d\,(d+2)}\,(\delta x)^{D_W}$, which allows a time step~$\displaystyle\frac{2\,d}{d+2}$ times larger than the one of the finite difference scheme. Indeed, since~$0<\lambda\leq d+2$ for any eigenvalue~$\lambda$ of~$\tilde{L}_m$, one has~$1-h\,\displaystyle\frac{d}{2}\,(d+2)^{m+1}\leq\gamma_i<1$.
	\end{corollary}
	
	\vskip 1cm
	
	\paragraph{Convergence}

	\begin{definition}
		\noindent
		\begin{enumerate}
			\item[\emph{i}.]  The scheme will be said to be \emph{convergent} if:
			{

			\begin{equation}  \displaystyle \lim_{h \to 0^+, \, m \to+\infty} \left  \| \left(u_j(t_n) -u^n_j \right)_{0\leq n\leq N,\, 1\leq j \leq d^m} \right  \|=0
				\, \cdot 	\end{equation}
			}
			\noindent where the~$u_j(t_n)$ denote the mean values of the exact solution, and the~$u^n_j$ the values given by the scheme.
			\item[\emph{ii}.] The scheme will be said to be \emph{conditionally convergent} if there exist two strictly positive constants~$\beta$ and $C$ such that, with~$\delta x=2^{-m}$:

		\begin{equation}\displaystyle \lim_{h\leq C\,(\delta x)^{\beta}, \, m \to+\infty} \left  \| \left(u_j(t_n) -u^n_j \right)_{0\leq n\leq N,\, 1\leq j \leq d^m} \right  \|=0
			\, \cdot 	\end{equation}

		\end{enumerate}
	\end{definition}

	\vskip 1cm

	\begin{theorem}{\ }\\

		{	\noindent Under assumptions~(H1) and~(H2), and under the CFL condition~(\ref{eq:CFLFVM}), the explicit finite volume scheme is convergent for the norm~$\| \cdot  \|_{2,\infty}$ defined by:

			\begin{equation}	\left \| \left( v^n_j \right)_{0\leq n\leq N, 1\leq j \leq d^m} \right  \|_{2,\infty}=\displaystyle \max_{0\leq n\leq N} \left( \displaystyle \frac{1}{d^m}\, \sum_{1\leq j\leq d^m}  \left  |v^n_j \right |^2\right)^{\frac{1}{2}} 	\end{equation}

			\noindent More precisely:

			$$\left\| \left(u_j(t_n) -u^n_j \right)\right\|_{2,\infty}=\mathcal{O}(h)+\mathcal{O}(2^{-m\alpha})=\mathcal{O}(2^{-m\alpha}) \, \cdot$$
		}
	\end{theorem}

	\vskip 1cm

	\begin{proof}

		\noindent Let us set:

		\begin{equation} w^n_j  = u_j(t_n) - u^n_j, \quad 0\leq n \leq N,\, 1\leq j\leq d^m 	\end{equation}

		\noindent and:

	\begin{equation}\forall\, n\, \in\, \left \lbrace 0, \hdots,   N \right \rbrace\, :
		\quad W^n=
		\left(
		\begin{matrix}
		w^n_1\\
		\vdots\\
		w^n_{d^m}
		\end{matrix}
		\right)
		\quad , \quad  E^n=
		\left(
		\begin{matrix}
		\varepsilon^m_{n,1}\\
		\vdots\\
		\varepsilon^m_{n,d^m}
		\end{matrix}
		\right)
			\end{equation}

		\noindent Since~$u^0_j=u_j(0)$, one has:

		\begin{equation}W^0=0	\end{equation}

		\noindent and, by definition of the consistency error:

		\begin{equation}\forall\, n\, \in\, \left \lbrace 0, \hdots,   N-1 \right \rbrace\, :
		\quad
		W^{n+1} =A\, W^n+h\, E^n 	\end{equation}

		\noindent By induction, one gets:

		\begin{equation}\forall\, n\, \in\, \left \lbrace 1, \hdots,   N \right \rbrace\, :
		\quad
		W^{n}=A^n W^0+h\, \displaystyle\sum_{k=0}^{n-1} A^{n-1-k}\,  E^{k}=h\, \displaystyle\sum_{k=0}^{n-1} A^{n-1-k} \, E^{k}
			\end{equation}

		\noindent Due to the symmetry of the matrix~$A$, the Courant-Friedrichs-Lewy stability condition
		$$h\,(d+2)^{m}\leq \displaystyle\frac{2}{d^2}$$

		\noindent  yields~$\rho(A)\leq 1$, and thus, denoting by~$|\cdot|$ the Euclidean norm, for any~$n\, \in\, \left \lbrace 1, \hdots,   N \right \rbrace$:

	\begin{equation}	\begin{array}{ccl}
		|W^n|&\leq &h\, \displaystyle\sum_{k=0}^{n-1}\rho(A)^{n-1-k}\,  | E^{k} |  \\
		&\leq & h\, n\,  \displaystyle \max_{0\leq k \leq N-1} | E^{k} | \\
		&\leq &h\, N\,  \displaystyle \max_{0\leq k \leq N-1} | E^{k} |  \\
		&=& T\,  \displaystyle \max_{0\leq k \leq N-1} \left(\sum_{j=1}^{d^m} \, |\varepsilon^m_{k,j}|^2 \right)^{\frac{1}{2}}  \\
		\end{array} 		\end{equation}

		\noindent and thus:

		\begin{align}
		\displaystyle 	\max_{0\leq n\leq N} \left(  \displaystyle \frac{1}{d^{m}} \, \sum_{j=1}^{d^m}   |w^n_j|^2\right)^{\frac{1}{2}}&= \displaystyle \frac{1}{d^{\frac{m}{2}}}\, \max_{1\leq n\leq N}|\, W^n|\\
		&\leq \displaystyle \frac{1}{d^{\frac{m}{2}}}\,   T\left(\max_{0\leq k \leq N-1} \left(\displaystyle\sum_{j=1}^{d^m} |\varepsilon^m_{k,j}|^2 \right)^{\frac{1}{2}} \right)\\
		&\leq \displaystyle \frac{1}{d^{\frac{m}{2}}} \, T\left(d^{\frac{m}{2}}\, \max_{0\leq k \leq N-1,\, 1\leq j \leq d^m} |\varepsilon^m_{k,j}| \right)\\
		&=  T\, \displaystyle\max_{0\leq k \leq N-1,\, 1\leq j \leq d^m} |\varepsilon^m_{k,j}|  \\
		&=   \mathcal{O}(h) + \mathcal{O}(2^{-m\alpha})
		\end{align}

		\noindent Under the CFL condition,~$h=\mathcal{O}((d+2)^{-m})$, and, since~$\alpha<1$, one has~$(d+2)^{-m}\leq 2^{-m}\leq 2^{-m\alpha}$, which yields the result.

	\end{proof}

	\vskip 1cm

	\begin{remark}
		\noindent Since the~$d^m$ $m^{th}$-order cells form a partition of~$\mathfrak{SS}$ (up to a set of measure zero), and since each of them has measure~$d^{-m}$, one has, for any function~$v$ which is constant on each~$m^{th}$-order cell, of value~$v_j$ on~$C^j_m$:
		{

		\begin{equation} \left( \displaystyle \frac{1}{d^m}\, \displaystyle \sum_{1\leq j\leq d^m}   |v_j|^2\right) ^{\frac{1}{2}}
			=\left \|  v \right \| _{L^2(\mathfrak{SS},\mu)} \, \cdot 	\end{equation}

		}
	\end{remark}

	\vskip 1cm
	
\subsubsection{The implicit Euler Method}

		\noindent Let us consider the implicit Euler scheme, for any integer~$n$ belonging to~$\left \lbrace 1, \hdots, N \right \rbrace$, and any integer~$j$ belonging to~$\left \lbrace 1, \hdots, d^m \right \rbrace$:

	\begin{equation}
	\label{eq:FVimplicit}
	u^{n}_j= u^{n-1}_j + \displaystyle \frac{h}{\mu(C^j_m)}\, \displaystyle \frac{1}{r^m}\frac{d}{2} \left(\sum_{l \underset{m}{\sim}j} \left(u^n_{l} - u^n_j\right)-2\,\varepsilon_j\,u^n_j\right) 	\, \cdot\end{equation}

		\noindent where, as in the explicit case,~$l\underset{m}{\sim}j$ means that the cells~$C^j_m$ and~$C^l_m$ are neighbors, and where~$\varepsilon_j=1$ if~$C^j_m$ contains a point of~$V_0$, and~$\varepsilon_j=0$ otherwise.\\

		\noindent It satisfies the recurrence relation:

	\begin{equation} \tilde{A} \, U(n) = U(n-1) \quad , \quad 1\leq n\leq N	\end{equation}

		\noindent where:

	\begin{equation}
		\tilde{A}=I_{d^m} + h\,\tilde{\Delta}_{m}=I_{d^m} + h\,\frac{d}{2}\,(d+2)^m\,\tilde{L}_m
			\end{equation}

		\noindent  and where~$I_{d^m}$ denotes the~$   d^m  \times   d^m $ identity matrix, and~$\tilde{L}_{m}$ the~$   d^m  \times   d^m $ matrix of the graph Laplacian of~$\mathcal{SS}_m$, with Dirichlet boundary terms, introduced in the explicit case.\\

		\noindent The scheme is well defined. Indeed, the matrix~$\tilde{A}$ is symmetric, and, thanks to proposition~\ref{prop:SpectralBoundFVM}, its eigenvalues~$\tilde{\gamma}_i$,~$1\leq i\leq d^m$, satisfy:

		$$1<\tilde{\gamma}_i\leq 1+h\,d^2\,(d+2)^m \, \cdot$$

		\noindent The matrix~$\tilde{A}$ is thus invertible, and:

		\begin{equation}
		\label{eq:ImplicitSpectrumFVM}
		\rho\left(\tilde{A}^{-1}\right)=\max_{1\leq i\leq d^m}\,\left|\tilde{\gamma}_i^{-1}\right|<1 \, \cdot
		\end{equation}

		\vskip 1cm

		\paragraph{Consistency, stability and convergence}$\,$
		\vskip 0.5cm

		\noindent \textbf{\emph{i}. Consistency}\\

		\noindent The consistency error of the implicit Euler scheme is defined, for~$1\leq n \leq N$ and~$1\leq j \leq d^m$, by:

		$$\varepsilon^m_{n,j}=\frac{u_j(t_{n})-u_j(t_{n-1})}{h}-\frac{d}{2}\,(d+2)^m\left(\sum_{l \underset{m}{\sim}j}\left(u_l(t_n)-u_j(t_n)\right)-2\,\varepsilon_j\,u_j(t_n)\right)$$

		\noindent where the~$u_j(t)$ denote the mean values of the exact solution over the cells~$C^j_m$. Under assumptions~(H1) and~(H2), the time discretization estimate (with a backward difference) and proposition~\ref{prop:SpaceConsistency} yield:

		\begin{equation}
		\varepsilon^m_{n,j} =\mathcal{O}(h) + \mathcal{O}\left(\displaystyle \frac{1}{2^{\alpha\, m}} \right)
		\end{equation}

		\noindent uniformly with respect to~$n$ and~$j$. One has:

		$$ \displaystyle \lim_{h \to 0^+, \, m \to+\infty} \,\max_{n,j}\, |\varepsilon^m_{n,j}|=0 $$

		\noindent which yields the consistency of the scheme.

		\vskip 1cm

		\noindent \textbf{\emph{ii}. Stability}\\

		\noindent The scheme writes:

		\begin{equation}
		U(n)=\tilde{A}^{-1}\,U(n-1)=\tilde{A}^{-n}\,U(0) \quad , \quad 1\leq n\leq N \, \cdot
		\end{equation}

		\noindent Since~$\tilde{A}^{-1}$ is symmetric, the spectral bound~(\ref{eq:ImplicitSpectrumFVM}) yields, for any~$h>0$ and any~$m\,\in\,\N^\star$:

		$$\rho\left(\tilde{A}^{-n}\right)=\rho\left(\tilde{A}^{-1}\right)^n\leq 1 \quad , \quad 1\leq n\leq N$$

		\noindent and thus, denoting by~$|\cdot|$ the Euclidean norm:

		\begin{equation} \left|U(n)\right| \leq  \left|U(0)\right| \quad , \quad 1\leq n\leq N\, \cdot \end{equation}

		\noindent The implicit scheme is thus unconditionally stable, in the sense of the definition given in the explicit case, with~$C=1$: no condition linking~$h$ and~$\delta x=2^{-m}$ is required.\\

		\noindent One may also note that~$\tilde{A}$ is a strictly diagonally dominant M-matrix. Thus,~$\tilde{A}^{-1}$ has nonnegative coefficients, and:

	\begin{equation} \parallel \tilde{A}^{-1} \parallel_{\infty} \leq 1\quad \Longrightarrow \quad
		\parallel \tilde{A}^{-n} \parallel_{\infty} \leq 1 	\end{equation}

		\noindent i.e.~$\parallel U(n)\parallel_{\infty} \leq \parallel U(0) \parallel_{\infty}$: the implicit scheme also satisfies the discrete maximum principle.

		\vskip 1cm
		\newpage

		\noindent \textbf{\emph{iii}. Convergence}\\

		\begin{theorem}
			\noindent Under assumptions~(H1) and~(H2), the implicit Euler scheme is convergent for the norm~$\parallel \cdot \parallel_{2,\infty}$, without any condition linking~$h$ and~$m$:

			$$\left\|\left(u_j(t_n)-u^n_j\right)_{0\leq n\leq N,\,1\leq j\leq d^m}\right\|_{2,\infty}=\mathcal{O}(h)+\mathcal{O}\left(2^{-m\alpha}\right) \, \cdot$$

		\end{theorem}

		\vskip 0.5cm

		\begin{proof}

			\noindent Let:

		\begin{equation} w^n_j  = u_j(t_n) - u^n_j, \quad 0\leq n \leq N,\, 1\leq j\leq d^m \, \cdot \end{equation}

			\noindent We set:

		\begin{equation} W^n=
			\left(
			\begin{matrix}
			w^n_1\\
			\vdots\\
			w^n_{d^m}
			\end{matrix}
			\right)
			\quad , \quad
			E^n=
			\left(
			\begin{matrix}
			\varepsilon^m_{n,1}\\
			\vdots\\
			\varepsilon^m_{n,d^m}
			\end{matrix}
			\right)
		\end{equation}

			\noindent Since~$u^0_j=u_j(0)$, one has~$W^0=0$. By definition of the consistency error, the exact mean values satisfy~$\tilde{A}\,\left(u_j(t_n)\right)_j=\left(u_j(t_{n-1})\right)_j+h\,E^n$, and thus, for~$1\leq n \leq N $:

			\begin{equation}
			W^{n}=\tilde{A}^{-1}\left( W^{n-1}+h\, E^n\right) \, \cdot
			\end{equation}

			\noindent Due to the stability of the scheme,~$\rho(\tilde{A}^{-1})\leq 1$, and we have, for~$n=1,\hdots,N$:

			\begin{align}
			|W^n|&\leq |W^{n-1}|+h\,|E^n|\\
			&\leq h\, \sum_{k=1}^{n}|E^k|\\
			&\leq h\,N\, \max_{1\leq k \leq N} | E^{k} | \\
			&= T\, \max_{1\leq k \leq N} \left(\sum_{j=1}^{d^m} |\varepsilon^m_{k,j}|^2 \right)^{\frac{1}{2}}
			\end{align}

			\noindent Thus:

			\begin{align}
			\displaystyle 	\max_{0\leq n\leq N} \left(  \frac{1}{d^m}\, \displaystyle \sum_{j=1}^{d^m}   |w^n_j|^2\right)^{\frac{1}{2}}&= \displaystyle \frac{1}{d^{\frac{m}{2}}}\, \max_{1\leq n\leq N}|W^n|\\
			& \leq  \displaystyle \frac{1}{d^{\frac{m }{2}}} \, T\, \max_{1\leq k \leq N} \, \left(\sum_{j=1}^{d^m} |\varepsilon^m_{k,j}|^2 \right)^{\frac{1}{2}} \\
			&\leq  \displaystyle \frac{1}{d^{\frac{m }{2}}}  \, T\, d^{\frac{m}{2}}\, \max_{1\leq k \leq N,\, 1\leq j \leq d^m} \, |\varepsilon^m_{k,j}|  \\
			&=   T \, \max_{1\leq k \leq N,\, 1\leq j \leq d^m} \, |\varepsilon^m_{k,j}|  \\
			&=   \mathcal{O}(h) + \mathcal{O}\left( \displaystyle \frac{1}{2^{\alpha\, m}}\right)
			\end{align}

			\noindent which goes to zero when~$h\to 0$ and~$m\to+\infty$, independently. The scheme is thus convergent.
		\end{proof}

		\vskip 1cm

		\begin{remark}
			\noindent With~$\delta x=2^{-m}$, the error of the implicit scheme writes~$\mathcal{O}(\delta t)+\mathcal{O}\left((\delta x)^{\frac{D_W-D_H}{2}}\right)$, where~$D_H=\displaystyle\frac{\ln d}{\ln 2}$ and~$D_W=\displaystyle\frac{\ln(d+2)}{\ln 2}$ denote the self-similar and walk dimensions of~$\mathfrak{SS}$. The two contributions are balanced for~$\delta t\asymp (\delta x)^{\frac{D_W-D_H}{2}}$, which is much larger than the time step~$\delta t\lesssim (\delta x)^{D_W}$ required by the explicit scheme. For a given accuracy, the implicit scheme thus requires far fewer time steps, at the cost of the resolution of a sparse, symmetric positive definite linear system at each step.
		\end{remark}

	\vskip 1cm

\subsubsection{Numerical application}

	\paragraph{The graph Laplacian sequence}$\,$\\

	\noindent In the sequel, we describe our recursive algorithm used to construct the matrices~$\tilde{L}_m$ related to the sequence of cell graphs~$\mathcal{SS}_m$, in the case of the Sierpi\'{n}ski Gasket and Tetrahedron. The vertices of~$\mathcal{SS}_m$ are the~$d^m$ cells~$f_{\cal W}(\mathfrak{SS})$,~$|{\cal W}|=m$, numbered according to the lexicographic order of the words~${\cal W}$. Contrary to the case of the finite difference method, no vertex is shared by two copies: the fusion of two copies consists in adding an edge between two cells.\\

	\noindent For any~$L\,\in\,\N$, and any~$i\,\in\,\{1,\hdots,d\}$, we denote by~$I_i(L)$ the label of the cell~$f_i\circ\hdots\circ f_i(\mathfrak{SS})$ (word~$(i,\hdots,i)$ of length~$L$), i.e. of the corner cell of~$\mathcal{SS}_L$ which contains the point~$P_i$:

	\begin{equation}
	I_i(0)=1 \quad , \quad \forall\,L\,\in\,\N^\star\,:\quad I_i(L)=I_i(L-1)+(i-1)\,d^{L-1} \quad , \quad \text{i.e.}\quad I_i(L)=1+(i-1)\,\frac{d^L-1}{d-1} \, \cdot
	\end{equation}

	\subparagraph{\emph{i}. The 2-Simplex (Triangle).}$\,$

	\vskip 1cm

	For any~\mbox{$m\,\in\,\N^\star$}, we will denote by $\text{Corner}_{k}$,~\mbox{$1\leq k \leq 3$}, the three corner cells of a~$m^{th}$-order triangle, i.e. the three cells which contain its vertices~(see figure~\ref{CornersSierpinskiGasket}). The~$(m+1)^{th}$-order triangle is then constructed by connecting three copies~$T_k$,~$k=1,\,2,\, 3$, of the~$m^{th}$-order triangle.\\

	\begin{center}
		\includegraphics[scale=0.5]{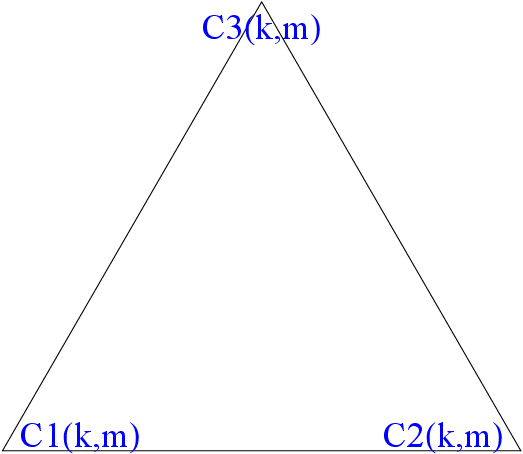}
		\captionof{figure}{$m^{th}$-order cell.}
		\label{CornersSierpinskiGasket}
	\end{center}

	\noindent The initial triangle is labelled such that $\text{Corner}_1 \sim 1$, $\text{Corner}_2 \sim 2$ and $\text{Corner}_3 \sim 3$~(see figure~\ref{CornersSierpinskiGasket}).

	\begin{figure}[!htb]
		\begin{center}
			\minipage{0.32\textwidth}%
			\includegraphics[width=\linewidth]{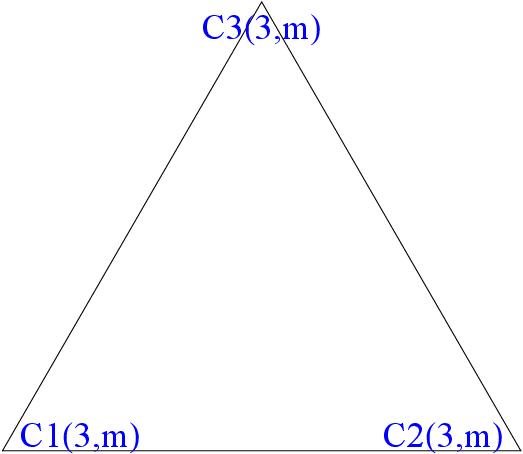}
			\caption{The third copy $T_3$}
			\label{figFVSG3}
			\endminipage
		\end{center}
		\minipage{0.32\textwidth}
		\includegraphics[width=\linewidth]{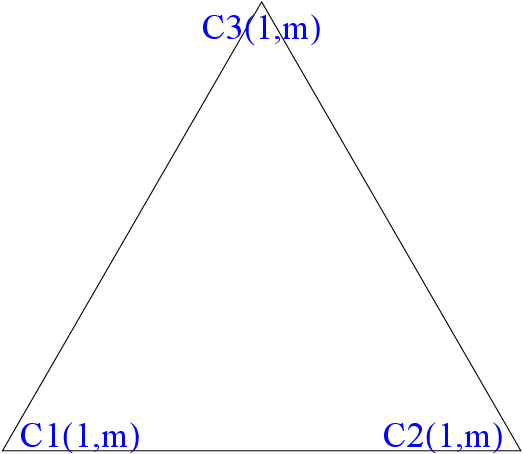}
		\caption{The first copy $T_1$}
		\label{figFVSG1}
		\endminipage\hfill
		\minipage{0.32\textwidth}
		\includegraphics[width=\linewidth]{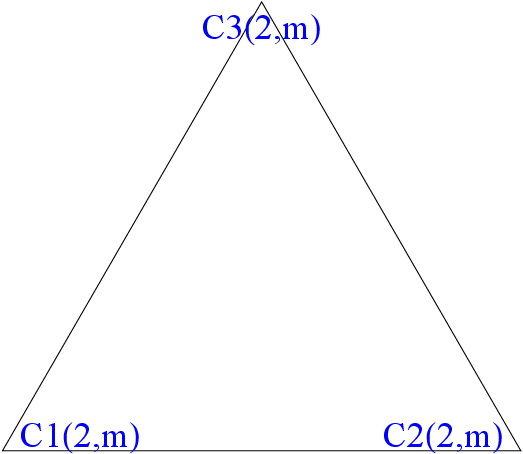}
		\caption{The second copy $T_2$}
		\label{figFVSG2}
		\endminipage\hfill
	\end{figure}

	\vskip 1cm

	\noindent \textbf{Iterative process:}\\

	\begin{enumerate}

		\item[\emph{i}.] The starting point is the cell graph~$\mathcal{SS}_1$, i.e. the complete graph on the three cells of order~$1$. The corresponding matrix is given by:

			\begin{equation} A_1 = \left(
		\begin{matrix}
		2 & -1 & -1 \\
		-1 & 2 & -1 \\
		-1 & -1 & 2 \\
		\end{matrix}
		\right) \end{equation}

		\item[\emph{ii}.] For any integer~$m\geq 2$, we first build the block diagonal matrix:

		$$B_m = \left(
		\begin{matrix}
		A_{m-1} & 0 & 0 \\
		0 & A_{m-1} & 0 \\
		0 & 0 & A_{m-1} \\
		\end{matrix}
		\right)$$

		\noindent the~$k^{th}$ block of which corresponds to the copy~$T_k$, whose cells are labelled from~$(k-1)\,3^{m-1}+1$ to~$k\,3^{m-1}$.

		\item[\emph{iii}.]  The fusion requires the following connections:

			\begin{equation} \left \lbrace
		\begin{array}{ccc}
		\text{Corner}_2(1,m)&\sim& \text{Corner}_1(2,m)\\
		\text{Corner}_3(1,m)& \sim & \text{Corner}_1(3,m)\\
		\text{Corner}_3(2,m) &\sim &\text{Corner}_2(3,m)\\
		\end{array}\right. \,\cdot \end{equation}

		\noindent where the labels of the corner cells of the copy~$T_k$ are given by:

			\begin{equation} \left \lbrace
		\begin{array}{ccc}
		\text{Corner}_1(k,m)&=&1+(k-1)\,3^{m-1}\\
		\text{Corner}_2(k,m)&=&I_2(m-1)+(k-1)\,3^{m-1}\\
		\text{Corner}_3(k,m)&=&k \,3^{m-1}\\
		\end{array}\right. \quad , \quad I_2(m-1)=\frac{3^{m-1}+1}{2}\,\cdot \end{equation}

		\item[\emph{iv}.] The connection matrix~(we refer to~\cite{UtaFreiberg2004}) is:

	\begin{equation} C_m = \left(
		\begin{matrix}
		\text{Corner}_2(1,m) & \text{Corner}_3(1,m) & \text{Corner}_3(2,m) \\
		\text{Corner}_1(2,m) & \text{Corner}_1(3,m) & \text{Corner}_2(3,m) \\
		\end{matrix}
		\right)\end{equation}

		\noindent For each column~$j$ of~$C_m$, the cells~$a=C_m(1,j)$ and~$b=C_m(2,j)$ are linked by a new edge: the matrix~$A_m$ is obtained from~$B_m$ by setting

	\begin{equation}  A_{m,\,ab}=A_{m,\,ba}=-1 \quad , \quad A_{m,\,aa}=B_{m,\,aa}+1=3 \quad , \quad A_{m,\,bb}=B_{m,\,bb}+1=3 \, \cdot \end{equation}

		\item[\emph{v}.] Finally, the matrix~$\tilde{L}_m$ of the scheme is obtained by adding the Dirichlet terms: the diagonal coefficients of the three corner cells~$I_1(m)=1$,~$I_2(m)$ and~$I_3(m)=3^m$ of~$\mathcal{SS}_m$ are increased by~$2$ (they thus become equal to~$4$).

	\end{enumerate}

	\vskip 1cm

	\subparagraph{\emph{ii}. The~$3$-Simplex (Tetrahedron).}$\,$

	\vskip 1cm

	For any~\mbox{$m\,\in\,\N^\star$}, we will denote by $\text{Corner}_{k}$,~\mbox{$1\leq k \leq 4$}, the four corner cells of a~$m^{th}$-order tetrahedron~(see figure~\ref{CornersSierpinskiTetrahedron}). The~$(m+1)^{th}$-order tetrahedron is constructed by connecting four copies~$T_k$,~$k=1,\hdots,4$, of the~$m^{th}$-order tetrahedron.\\

	\begin{center}
		\includegraphics[scale=0.7]{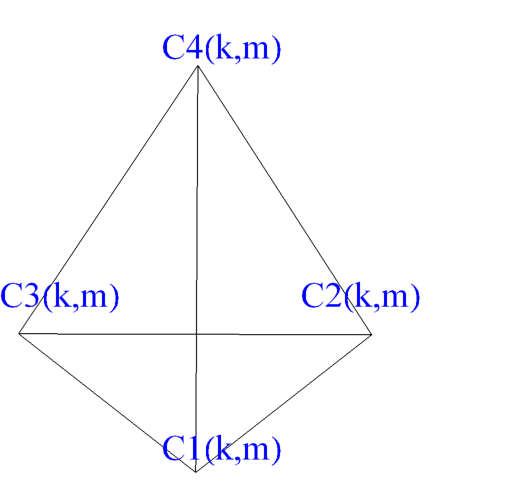}
		\captionof{figure}{$m^{th}$-order cell.}
		\label{CornersSierpinskiTetrahedron}
	\end{center}

	\begin{figure}[!htb]

		\begin{center}
			\minipage{0.32\textwidth}%
			\includegraphics[width=\linewidth]{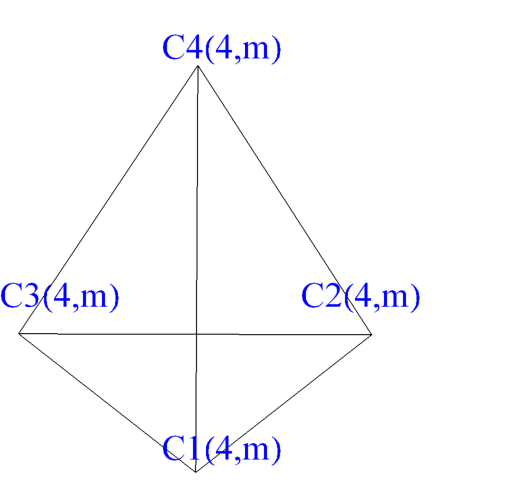}
			\caption{The fourth copy $T_4$.}
			\label{figFVST4}
			\endminipage
		\end{center}
		\minipage{0.32\textwidth}
		\includegraphics[width=\linewidth]{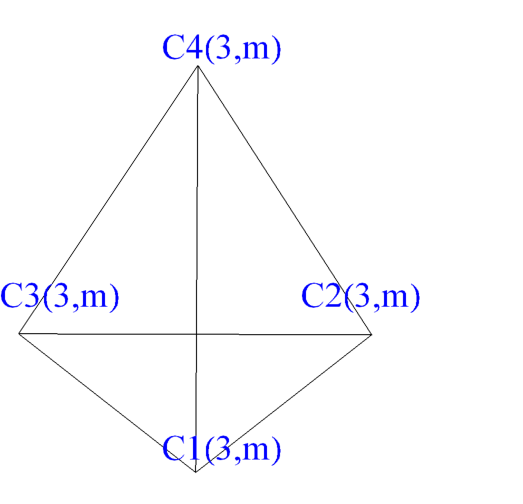}
		\caption{The third copy $T_3$.}
		\label{figFVST3}
		\endminipage\hfill
		\minipage{0.32\textwidth}
		\includegraphics[width=\linewidth]{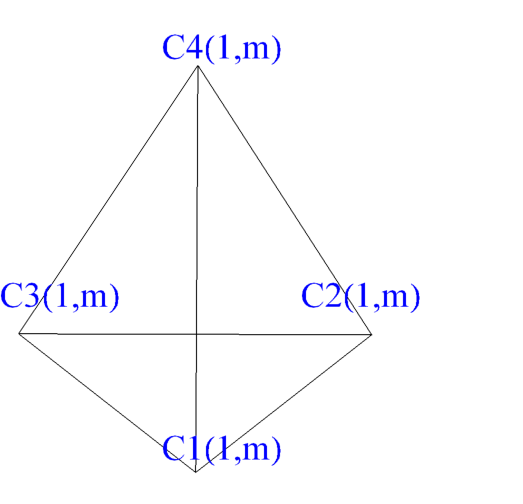}
		\caption{The first copy $T_1$.}
		\label{figFVST1}
		\endminipage\hfill
		\minipage{0.32\textwidth}
		\includegraphics[width=\linewidth]{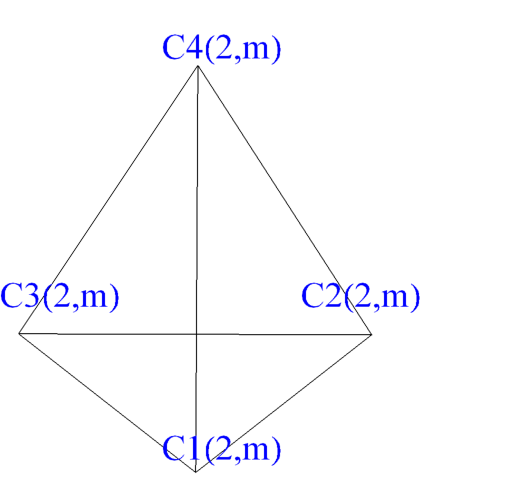}
		\caption{The second copy $T_2$.}
		\label{figFVST2}
		\endminipage\hfill

	\end{figure}








	\noindent For the fusion, one has to make the following connections:

\begin{equation}\left \lbrace
	\begin{array}{ccc}
	\text{Corner}_2(1,m)&\sim& \text{Corner}_1(2,m)\\
	\text{Corner}_3(1,m)& \sim & \text{Corner}_1(3,m)\\
	\text{Corner}_4(1,m) &\sim & \text{Corner}_1(4,m)\\
	\text{Corner}_3(2,m) &\sim &\text{Corner}_2(3,m)\\
	\text{Corner}_4(2,m) &\sim &\text{Corner}_2(4,m)\\
	\text{Corner}_4(3,m) &\sim &\text{Corner}_3(4,m)
	\end{array}\right. \,\cdot \end{equation}

	\vskip 1cm

	\begin{pte}\textbf{Labels of the corner cells of the copy~$T_k$,~$m\geq 2$}\\

		\noindent One has:
	\begin{equation} \left \lbrace
		\begin{array}{ccc}
		\text{Corner}_1(k,m)&=&1+(k-1)\,4^{m-1}\\
		\text{Corner}_2(k,m)&=&I_2(m-1)+(k-1)\,4^{m-1}\\
		\text{Corner}_3(k,m)&=&I_3(m-1)+(k-1) \,4^{m-1}\\
		\text{Corner}_4(k,m)&=&k \,4^{m-1}\\
		\end{array}\right. \,\cdot \end{equation}

		\noindent where:

		$$I_2(m-1)=1+\frac{4^{m-1}-1}{3}\quad , \quad I_3(m-1)=1+2\,\frac{4^{m-1}-1}{3}\, \cdot$$

		\noindent Equivalently:~$I_2(1)=2$,~$I_3(1)=3$, and, for any~$L\geq 2$,~$I_2(L)=I_2(L-1)+4^{L-1}$,~$I_3(L)=I_3(L-1)+2\times 4^{L-1}$.

	\end{pte}

	\vskip 1cm

	\noindent \textbf{Iterative process:}\\

	\begin{enumerate}

		\item[\emph{i}.] The starting point is the cell graph~$\mathcal{SS}_1$, i.e. the complete graph on the four cells of order~$1$. One has:

	\begin{equation}A_1 = \left(
		\begin{matrix}
		3 & -1 & -1 & -1\\
		-1 & 3 & -1 & -1\\
		-1 & -1 & 3 & -1 \\
		-1 & -1 & -1 & 3 \\
		\end{matrix}
		\right) \, \cdot \end{equation}
		\vskip 0.5cm

		\item[\emph{ii}.] For any integer~$m\geq 2$, we build the block diagonal matrix~$B_m$, made of four copies of~$A_{m-1}$, the~$k^{th}$ block corresponding to the copy~$T_k$.

		\vskip 0.5cm

		\item[\emph{iii}.] The connection matrix is:

		\noindent $\forall\, m\,\geq\,2\, :$\\
		\footnotesize
		{	$$C_m = \left(
			\begin{matrix}
			\text{Corner}_2(1,m) & \text{Corner}_3(1,m) & \text{Corner}_3(2,m) & \text{Corner}_4(1,m) & \text{Corner}_4(2,m) & \text{Corner}_4(3,m)\\
			\text{Corner}_1(2,m) & \text{Corner}_1(3,m) & \text{Corner}_2(3,m) & \text{Corner}_1(4,m) & \text{Corner}_2(4,m) & \text{Corner}_3(4,m)\\
			\end{matrix}
			\right)$$
			\normalsize

			\vskip 0.5cm

			\noindent For each column~$j$ of~$C_m$, the cells~$a=C_m(1,j)$ and~$b=C_m(2,j)$ are linked by a new edge:
			$$A_{m,\,ab}=A_{m,\,ba}=-1 \quad , \quad A_{m,\,aa}=B_{m,\,aa}+1=4 \quad , \quad A_{m,\,bb}=B_{m,\,bb}+1=4 \, \cdot$$}

		\item[\emph{iv}.] Finally, the matrix~$\tilde{L}_m$ of the scheme is obtained by increasing by~$2$ the diagonal coefficients of the four corner cells~$I_1(m)$,~$I_2(m)$,~$I_3(m)$,~$I_4(m)$ of~$\mathcal{SS}_m$ (Dirichlet terms).

	\end{enumerate}

	\vskip 0.5cm

	\noindent We have checked numerically, for~$d\in\{3,4\}$ and~$1\leq m\leq 4$, that this algorithm yields exactly the matrix of the graph Laplacian of~$\mathcal{SS}_m$, with the Dirichlet terms.

	\vskip 1cm

	\paragraph{Numerical results} $\, $
	\vskip 0.5cm

	 Our heat transfer simulation consists in a propagation scenario, where the initial condition is a harmonic spline~$g$ of level~$m$, which takes the value~$1$ on a vertex~$X$ of~$V_m\setminus V_0$, and~$0$ on the other vertices of~$V_m$. Its support is the union of the two~$m^{th}$-order cells which contain~$X$, and the initial values of the scheme are its mean values~$u^0_j$, which are equal to~$\displaystyle\frac{1}{d}$ on these two cells, and to~$0$ elsewhere.\\

	  Each point represents an~$m^{th}$-order cell. The colour scale represents the value of the temperature, from red (high values) to blue (low values).\\

	\noindent In both cases, the time step satisfies the stability condition of the corollary on the sharp stability condition, i.e.~$h\,(d+2)^{m+1}\leq\displaystyle\frac{4}{d}$.\\

	\subparagraph{\emph{i}. The~$2$-Simplex (Triangle)}$\,$\\

	\noindent In the following (see figures~\ref{figFVSG0} to~\ref{figFVSG500}), we give numerical results in the case where:

	$$m=6\quad , \quad T=1 \quad , \quad N=70.5\times 10^3\, ,$$

	\noindent i.e.~$h\approx 1.42\times 10^{-5}$, while the stability condition requires~$h\leq\displaystyle\frac{4}{3\times 5^7}\approx 1.71\times 10^{-5}$.

	\begin{center}
		\includegraphics[scale=1.5]{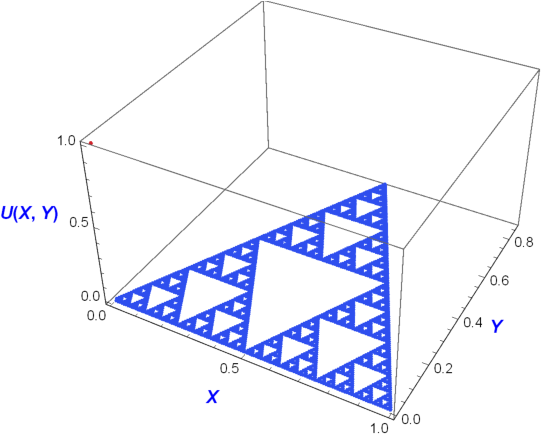}
		\captionof{figure}{The numerical solution for the initial condition.}
		\label{figFVSG0}
	\end{center}

	\begin{center}
		\includegraphics[scale=1.5]{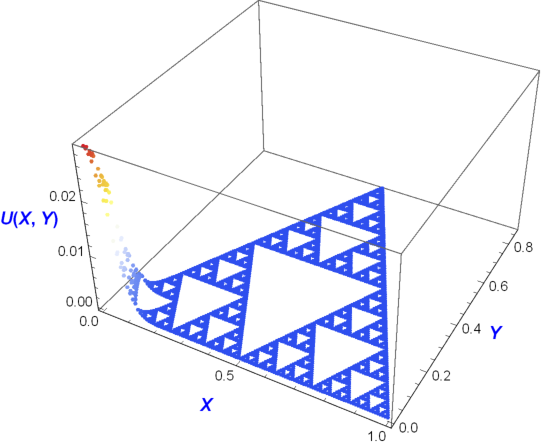}
		\captionof{figure}{The numerical solution, in the case~$n=100$.}
		\label{figFVSG100}
	\end{center}

	\begin{center}
		\includegraphics[scale=1.5]{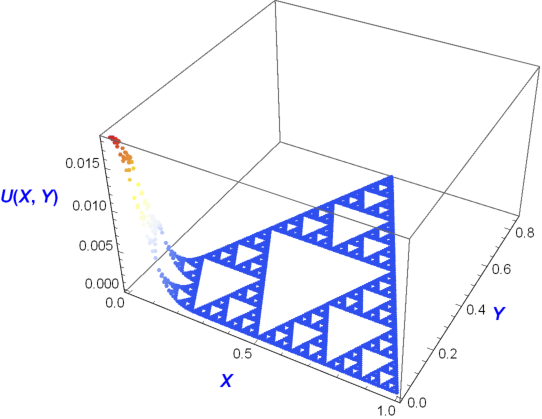}
		\captionof{figure}{The numerical solution, in the case~$n=200$.}
		\label{figFVSG200}
	\end{center}

	\begin{center}
		\includegraphics[scale=1.5]{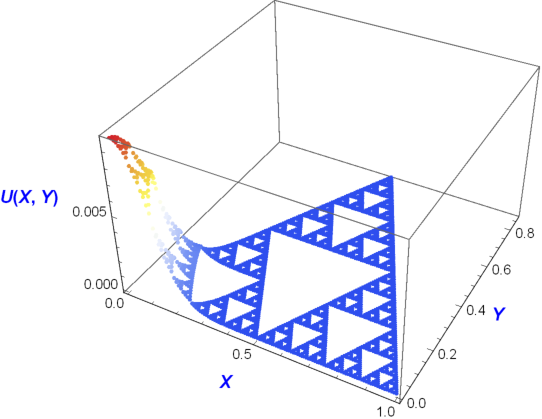}
		\captionof{figure}{The numerical solution, in the case~$n=500$.}
		\label{figFVSG500}
	\end{center}

	\newpage
	\subparagraph{\emph{ii}. The 3-Simplex (Tetrahedron)}$\,$\\

	\noindent In the following (see figures~\ref{figFVST0} to~\ref{figFVST500}), we give a four dimensional representation (three space coordinates and the colour scale) of the numerical results in the case where:

	$$m=5\quad , \quad T=1 \quad , \quad N=60\times 10^3\, ,$$

	\noindent i.e.~$h\approx 1.67\times 10^{-5}$, while the stability condition requires~$h\leq\displaystyle\frac{4}{4\times 6^6}\approx 2.14\times 10^{-5}$.

	\begin{center}
		\minipage{0.5\textwidth}
		\includegraphics[width=\linewidth]{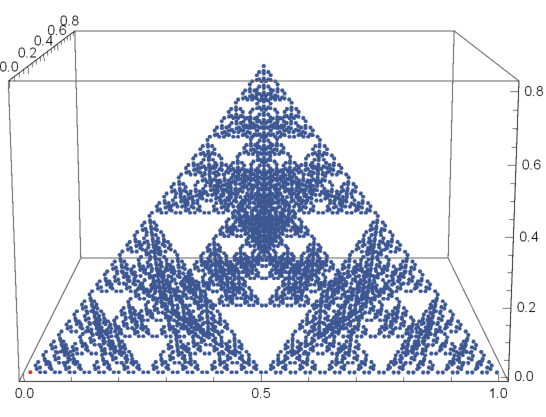}
		\endminipage\hfill
		\minipage{0.5\textwidth}
		\includegraphics[width=\linewidth]{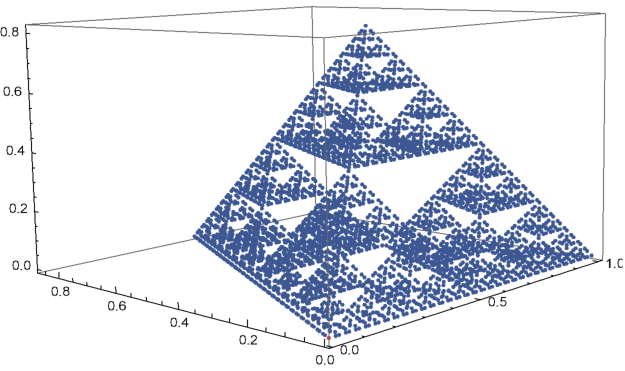}
		\endminipage\hfill
		\captionof{figure}{The numerical solution for the initial condition (two views).}
		\label{figFVST0}
	\end{center}

	\begin{center}
		\minipage{0.5\textwidth}
		\includegraphics[width=\linewidth]{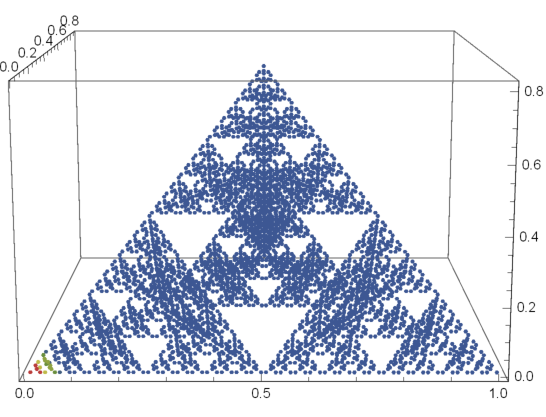}
		\endminipage\hfill
		\minipage{0.5\textwidth}
		\includegraphics[width=\linewidth]{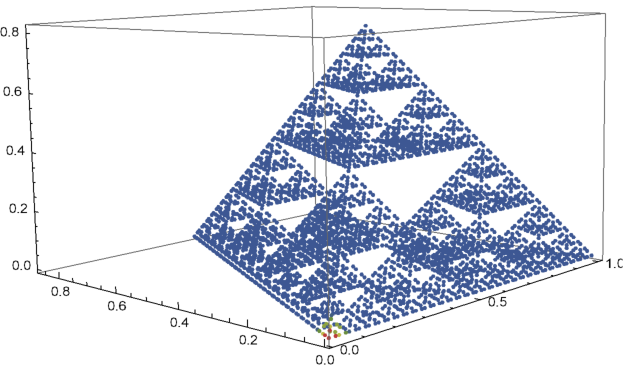}
		\endminipage\hfill
		\captionof{figure}{The numerical solution, in the case~$n=10$.}
		\label{figFVST10}
	\end{center}

	\begin{center}
		\minipage{0.5\textwidth}
		\includegraphics[width=\linewidth]{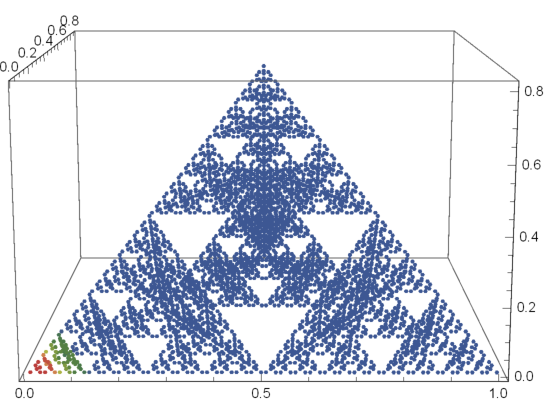}
		\endminipage\hfill
		\minipage{0.5\textwidth}
		\includegraphics[width=\linewidth]{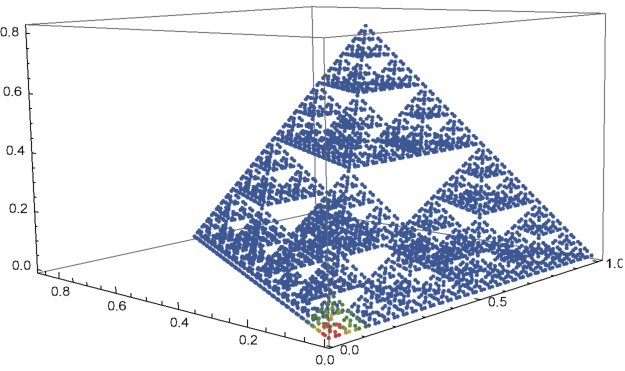}
		\endminipage\hfill
		\captionof{figure}{The numerical solution, in the case~$n=50$.}
		\label{figFVST50}
	\end{center}

	\begin{center}
		\minipage{0.5\textwidth}
		\includegraphics[width=\linewidth]{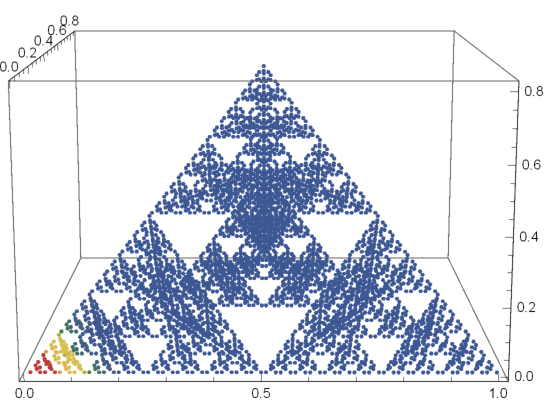}
		\endminipage\hfill
		\minipage{0.5\textwidth}
		\includegraphics[width=\linewidth]{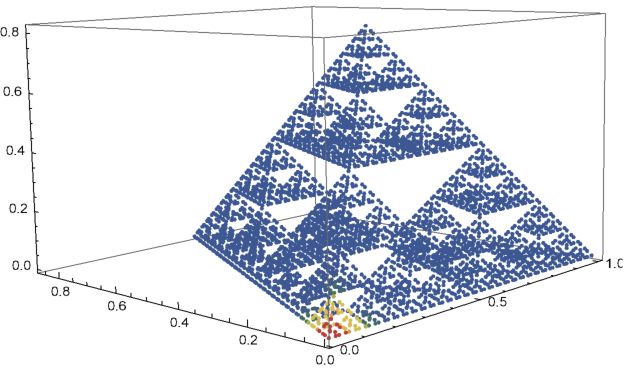}
		\endminipage\hfill
		\captionof{figure}{The numerical solution, in the case~$n=100$.}
		\label{figFVST100}
	\end{center}

	\begin{center}
		\minipage{0.5\textwidth}
		\includegraphics[width=\linewidth]{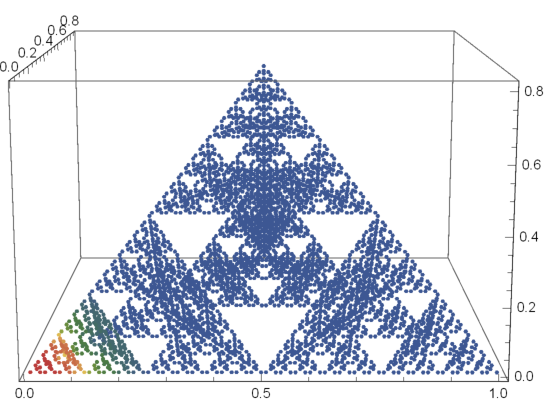}
		\endminipage\hfill
		\minipage{0.5\textwidth}
		\includegraphics[width=\linewidth]{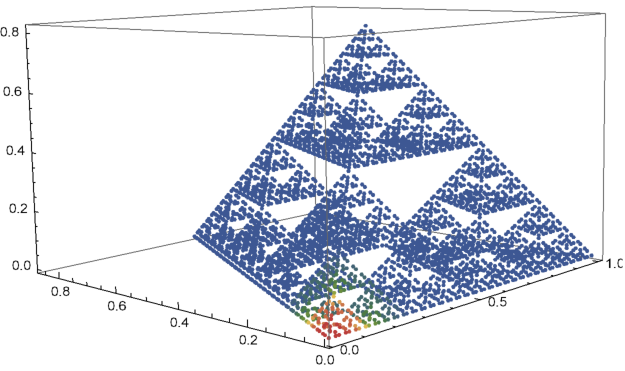}
		\endminipage\hfill
		\captionof{figure}{The numerical solution, in the case~$n=200$.}
		\label{figFVST200}
	\end{center}

	\begin{center}
		\minipage{0.5\textwidth}
		\includegraphics[width=\linewidth]{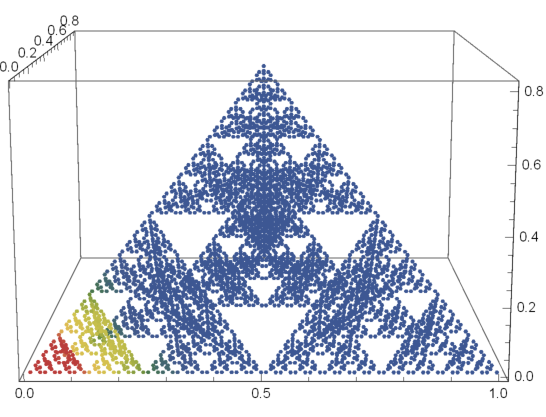}
		\endminipage\hfill
		\minipage{0.5\textwidth}
		\includegraphics[width=\linewidth]{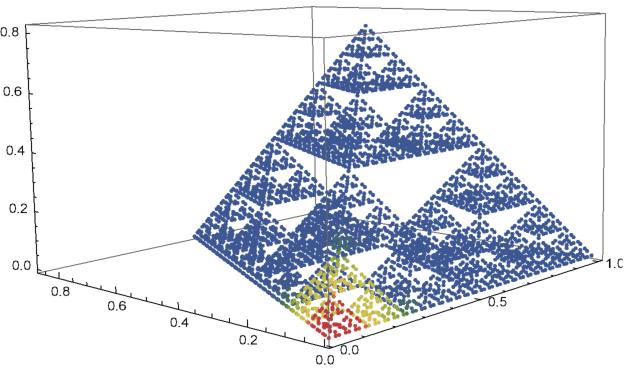}
		\endminipage\hfill
		\captionof{figure}{The numerical solution, in the case~$n=500$.}
		\label{figFVST500}
	\end{center}

	\vskip 1cm

	\subparagraph{\emph{iii}. Discussion}{$\,$}\\

	\noindent Let us consider the energy, i.e. the square of the~$L^2$ norm of the approximate solution, which is piecewise constant on the~$m^{th}$-order cells:

	$${E}(n)=\displaystyle  \int_{\mathfrak{SS}} u_h(t_n,X)^2 \, d\mu(X)=\frac{1}{d^m}\sum_{j=1}^{d^m}\left(u^n_j\right)^2 \, \cdot$$

	\noindent Plotted as a function of the iteration step~$n$ (see figure~\ref{figEnergy}), it exhibits an exponential decay. This is in accordance with the theory: expanding~$U(0)$ on an orthonormal basis of eigenvectors of~$\tilde{L}_m$, one gets, for large~$n$,~$E(n)\approx C\,\gamma_{\max}^{2n}$, where~$\gamma_{\max}=1-h\,\displaystyle\frac{d}{2}\,(d+2)^m\,\lambda_{\min}(\tilde{L}_m)$ is the largest eigenvalue of~$A$, i.e.~$E(n)\approx C\,e^{-2\,\tilde{\lambda}_1\,t_n}$, where~$\tilde{\lambda}_1=\displaystyle\frac{d}{2}\,(d+2)^m\,\lambda_{\min}(\tilde{L}_m)$ approximates the first Dirichlet eigenvalue of the Laplacian on~$\mathfrak{SS}$.


	\begin{center}
		\includegraphics[scale=1]{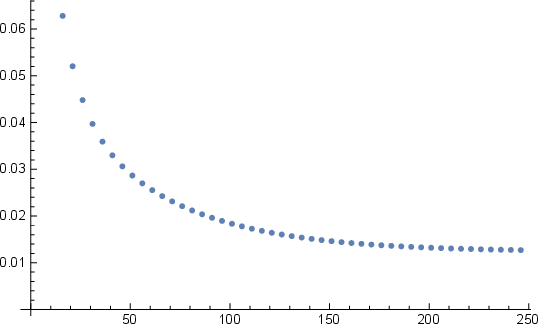}
		\captionof{figure}{The energy ${E}$ on the Sierpi\'{n}ski gasket as a function of $n$, for $m=3$.}
		\label{figEnergy}
	\end{center}

	\begin{center}
		\includegraphics[scale=1]{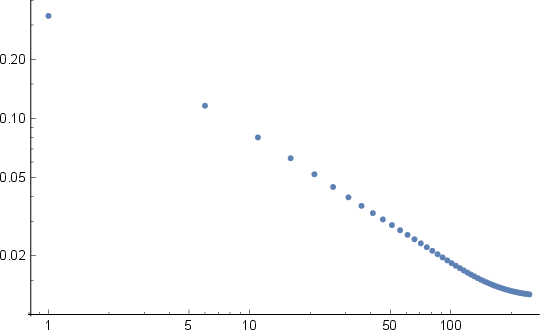}
		\captionof{figure}{The log-log plot of the energy ${E}$ on the Sierpi\'{n}ski gasket, for $m=3$.}
		\label{figEnergyLogLog}
	\end{center}

	\vskip 1cm

	By comparing the explicit finite difference method (FDM) and the finite volume method (FVM), one may deduce from the theoretical results the following similarities and differences:

	\begin{itemize}
		\item[\emph{i}.] The FDM is based on the sequence of graphs~\mbox{$\left (\mathfrak{SS}_m\right)_{m\in\N}$}, with~$\mathcal{N}_m-d=\displaystyle\frac{d^{m+1}-d}{2}$ unknowns, and the FVM is based on the sequence of cell graphs~\mbox{$\left (\mathcal{SS}_m\right)_{m\in\N}$}, with~$d^m$ unknowns. Both satisfy the spectral decimation property, with the same function~$\Phi$. Moreover, the spectrum of the FVM matrix is contained in the one of the FDM matrix (proposition~\ref{prop:SpectralInclusion}).
		\item[\emph{ii}.] The normalization constant~$\displaystyle\frac{d}{2}\,(d+2)^m$ is the same for both methods.
		\item[\emph{iii}.] The theoretical errors of both methods are of the same order,~$\mathcal{O}(h)+\mathcal{O}(2^{-m\alpha})$, under the same assumptions. In particular, both explicit schemes are of order~$1$ in time.
		\item[\emph{iv}.] The stability exponent is the same: in both cases, the explicit scheme is stable for~$h\lesssim(\delta x)^{D_W}$, where~$D_W$ is the walk dimension. However, the constant is better in the case of the FVM, which allows a time step~$\displaystyle\frac{2\,d}{d+2}$ times larger than the one of the FDM. The implicit schemes are unconditionally stable in both cases.
		\item[\emph{v}.] Finally, the numerical simulations show the same qualitative behavior in both approaches.
	\end{itemize}

\noindent The numerical method developed here is the continuation of our work~\cite{RianeDavidFDMS}. Our objective is to provide numerical approximations of the solutions of partial differential equations on fractals, since analytic solutions are difficult, or even impossible, to obtain in general. In the case of the heat equation, we have proved the convergence of the finite volume scheme, and described the behaviour of the approximate solution, which opens the way to the study of similar phenomena.

	\bibliographystyle{alpha}
	\bibliography{BibliographieClaire}
\end{document}